\documentclass[12pt,leqno]{article}
\textheight9in \def\DATE{\today}
\newtheorem{theorem}{Theorem}
\newtheorem{definition}[theorem]{Definition}

\newtheorem{example}[theorem]{Example}
\newtheorem{lemma}[theorem]{Lemma}

\newtheorem{proposition}[theorem]{Proposition}

\newtheorem{exercise}[theorem]{Exercise}

\catcode`\@=11

\def\ps@myheadings{\let\@mkboth\@gobbletwo
\def\@oddhead{\ifnum\count0=1 \hfill\else
\rightmark \hfil \rm\thepage\fi}%
\def\@oddfoot{\ifnum\count0=1 \hfill \rm 1 \hfill \else
\hfill\fi}
\def\@evenhead%
{\rm\leftmark\hfil\rm\thepage}%
\def\@evenfoot{}\def\sectionmark##1{}
\def\subsectionmark##1{}}

\def\@begintheorem#1#2{\it \trivlist \item[\hskip
 \labelsep{\bf #1\ #2.}]}
\def\@opargbegintheorem#1#2#3{\it \trivlist\item[\hskip%
 \labelsep{\bf #1\ #2.\ (#3)}]}
\def\@endtheorem{\endtrivlist}

\def\@listI{\leftmargin\leftmargini \parsep 1pt plus 2.5pt
 minus 1pt\topsep 5pt plus 2pt minus 3pt%
 \itemsep 0pt plus 2.5pt minus 1pt}
\let\@listi\@listI
\@listi

\def\@sect#1#2#3#4#5#6[#7]#8{\ifnum #2>\c@secnumdepth%
 \def \@svsec {}\else \refstepcounter {#1}\edef \@svsec%
 {\csname the#1\endcsname. \hskip .1em }\fi \@tempskipa%
 #5\relax \ifdim \@tempskipa >\z@ \begingroup #6\relax%
 \@hangfrom {\hskip #3\relax \@svsec }{\interlinepenalty%
 \@M #8.\par }\endgroup \csname #1mark\endcsname {#7}%
 \addcontentsline {toc}{#1}{\ifnum #2>\c@secnumdepth%
 \else \protect \numberline {\csname the#1\endcsname. }%
 \fi #7}\else \def \@svsechd {#6\hskip #3\@svsec #8.%
 \csname #1mark\endcsname {#7}\addcontentsline {toc}{#1}%
 {\ifnum #2>\c@secnumdepth \else \protect \numberline%
 {\csname the#1\endcsname. }\fi #7}}\fi \@xsect {#5}}

\def\section{\@startsection {section}{1}{\z@ }%
 {-3.5ex plus -1ex minus -.2ex}{2.3ex plus .2ex}{\bf }}

\def\thebibliography#1{%
 \section *{References.\@mkboth {REFERENCES}{REFERENCES}}%
 \list {[\arabic {enumi}]}{\settowidth \labelwidth {[#1]}%
 \leftmargin \labelwidth \advance \leftmargin \labelsep %
 \usecounter {enumi}} \def \newblock %
 {\hskip .11em plus .33em minus -.07em} \sloppy \clubpenalty 4000%
 \widowpenalty 4000 \sfcode`\.=1000\relax}

\def\@maketitle{%
 \newpage \null \vskip 2em
 \begin{center}
{\Large\bf \@title \par }
 \vskip 1.5em
 {\large \lineskip .5em
 \begin {tabular}[t]{c}\@author
 \end{tabular}\par}
 \end{center}
  \vskip .8em}

\def\abstract{%
\if@twocolumn \section *{Abstract}
 \else \small\quotation\noindent{\bf Abstract.}\fi}

\catcode`\@=13

\def\pravystromecek{
\put(19,5){\line(1,1){6}}
\put(25,11){\line(1,-1){6}}
\put(26.8,9.2){\line(-1,-1){4}}
\put(28.2,7.8){\line(-1,-1){2.8}}
\put(28.2,7.8){\line(0,-1){2.8}}
\put(25,9.5){\multiput(0,0)(.3,-.3){3}{\scriptsize $\cdot$}}
\put(29.3,6.5){\makebox(0,0)[cc]{\scriptsize $\ddots$}}
\put(28,8.5){\makebox(0,0)[lb]{\scriptsize $e_{i-1}$}}
}

\def\levystromecek{
\put(12,11){\line(1,-1){6}}
\put(12,11){\line(-1,-1){6}}
\put(13.8,9.2){\line(-1,-1){4.2}}
\put(15.2,7.8){\line(-1,-1){2.8}}
\put(11.8,9.5){\multiput(0,0)(.3,-.3){3}{\scriptsize $\cdot$}}
\put(16,6.5){\makebox(0,0)[cc]{\scriptsize $\ddots$}}
\put(15,8.5){\makebox(0,0)[lb]{\scriptsize $e'_i$}}
}

\def\doubless#1#2{{
\def\arraystretch{.5}
\begin{array}{c}
\mbox{\scriptsize $\scriptstyle #1$}
\\
\mbox{\scriptsize $\scriptstyle #2$}
\end{array}\def\arraystretch{1}
}}

\def\squeezedcdots{\hskip -.5mm \cdot \hskip -.5mm  \cdot \hskip -.5mm \cdot}

\def\jednadva{{
\unitlength=.4pt
\begin{picture}(24.00,20.00)(-2.00,0.00)
\bezier{20}(10.00,10.00)(15.00,5.00)(20.00,0.00)
\bezier{20}(10.00,10.00)(5.00,5.00)(0.00,0.00)
\put(10.00,20.00){\line(0,-1){10.00}}
\end{picture}}
}

\def\jednactyri{{
\unitlength=.05pt
\begin{picture}(176.00,160.00)(-8.00,0.00)
\put(80.00,100.00){\line(0,1){60.00}}
\bezier{20}(80.00,80.00)(100.00,30.00)(110.00,0.00)
\bezier{20}(80.00,80.00)(60.00,30.00)(50.00,0.00)
\bezier{20}(80.00,80.00)(120.00,40.00)(160.00,0.00)
\bezier{20}(80.00,80.00)(40.00,40.00)(0.00,0.00)
\put(80.00,100.00){\line(0,-1){20.00}}
\end{picture}}
}

\def\ctyrijedna{{
\unitlength=.05pt
\begin{picture}(176.00,160.00)(-8.00,-160.00)
\put(80.00,-100.00){\line(0,-1){60.00}}
\bezier{20}(80.00,-80.00)(100.00,-30.00)(110.00,0.00)
\bezier{20}(80.00,-80.00)(60.00,-30.00)(50.00,0.00)
\bezier{20}(80.00,-80.00)(120.00,-40.00)(160.00,0.00)
\bezier{20}(80.00,-80.00)(40.00,-40.00)(0.00,0.00)
\put(80.00,-80.00){\line(0,-1){40.00}}
\end{picture}}
}

\def\dvajedna{{
\unitlength=.4pt
\begin{picture}(24.00,20.00)(-2.00,0.00)
\put(10.00,10.00){\line(0,-1){10.00}}
\bezier{20}(10.00,10.00)(15.00,15.00)(20.00,20.00)
\bezier{20}(0.00,20.00)(5.00,15.00)(10.00,10.00)
\end{picture}}
}

\def\dvadva{{
\unitlength=.8pt
\begin{picture}(12.00,10.00)(-1.00,0.00)
\bezier{30}(0.00,0.00)(5.00,5.00)(10.00,10.00)
\bezier{30}(0.00,10.00)(5.00,5.00)(10.00,0.00)
\end{picture}}
}

\def\jednatri{{
\unitlength=.4pt
\begin{picture}(24.00,20.00)(-2.00,0.00)
\bezier{20}(10.00,10.00)(15.00,5.00)(20.00,0.00)
\bezier{20}(10.00,10.00)(5.00,5.00)(0.00,0.00)
\put(10.00,20.00){\line(0,-1){20.00}}
\end{picture}}
}

\def\trijedna{{
\unitlength=.4pt
\begin{picture}(24.00,20.00)(-2.00,-20.00)
\bezier{20}(10.00,-10.00)(15.00,-5.00)(20.00,0.00)
\bezier{20}(10.00,-10.00)(5.00,-5.00)(0.00,0.00)
\put(10.00,-20.00){\line(0,1){20.00}}
\end{picture}}
}

\def\dvatri{{
\unitlength=0.4pt
\begin{picture}(24.00,20.00)(-2.00,0.00)
\put(10.00,10.00){\line(0,-1){10.00}}
\bezier{30}(0.00,0.00)(10.00,10.00)(20.00,20.00)
\bezier{30}(0.00,20.00)(10.00,10.00)(20.00,0.00)
\end{picture}}
}

\def\tridva{{
\unitlength=.4pt
\begin{picture}(24.00,20.00)(-2.00,-20.00)
\put(10.00,-10.00){\line(0,1){10.00}}
\bezier{30}(0.00,0.00)(10.00,-10.00)(20.00,-20.00)
\bezier{30}(0.00,-20.00)(10.00,-10.00)(20.00,0.00)
\end{picture}}
}

\def\ZbbZbZb{{
\unitlength=.05pt
\begin{picture}(192.00,120.00)(-16.00,0.00)
\bezier{10}(30.00,30.00)(50.00,10.00)(60.00,0.00)
\bezier{20}(50.00,50.00)(70.00,30.00)(100.00,0.00)
\bezier{30}(80.00,80.00)(120.00,40.00)(160.00,0.00)
\bezier{30}(80.00,80.00)(40.00,40.00)(0.00,0.00)
\put(80.00,160.00){\line(0,-1){80.00}}
\end{picture}}
}

\def\bZbZbbZ{{
\unitlength=0.05pt
\begin{picture}(192.00,160.00)(-16.00,0.00)
\bezier{10}(100.00,0.00)(120.00,20.00)(130.00,30.00)
\bezier{20}(110.00,50.00)(90.00,30.00)(60.00,0.00)
\put(80.00,120.00){\line(0,1){40.00}}
\bezier{30}(80.00,80.00)(120.00,40.00)(160.00,0.00)
\bezier{30}(80.00,80.00)(40.00,40.00)(0.00,0.00)
\put(80.00,120.00){\line(0,-1){40.00}}
\end{picture}}
}

\def\gen#1#2{
\if #11
    \if #22 \jednadva \else \fi
\else
\fi
\if #12
    \if #22 \dvadva \else \fi
\else
\fi
\if #12
    \if #21 \dvajedna \else \fi
\else
\fi
\if #13
    \if #22 \tridva \else \fi
\else
\fi
\if #13
    \if #21 \trijedna \else \fi
\else
\fi
\if #12
    \if #23 \dvatri \else \fi
\else
\fi
\if #11
    \if #23 \jednatri \else \fi
\else
\fi
\if #11
    \if #24 \jednactyri \else \fi
\fi
\if #14
    \if #21 \ctyrijedna \else \fi
\fi
}

\def\bZbbZ{
{
\unitlength=.27pt
\begin{picture}(48.00,30.00)(-4,0.00)
\bezier{34}(20.00,20.00)(30.00,10.00)(40.00,0.00)
\bezier{34}(20.00,20.00)(10.00,10.00)(0.00,0.00)
\bezier{20}(30.00,10.00)(25.00,5.00)(20.00,0.00)
\put(20.00,30.00){\line(0,-1){10.00}}
\end{picture}}
}

\def\ZbbZb{{
\unitlength=.27pt
\begin{picture}(48.00,30.00)(-4,0.00)
\bezier{34}(20.00,20.00)(30.00,10.00)(40.00,0.00)
\bezier{34}(20.00,20.00)(10.00,10.00)(0.00,0.00)
\bezier{20}(10.00,10.00)(15.00,5.00)(20.00,0.00)
\put(20.00,30.00){\line(0,-1){10.00}}
\end{picture}}}

\def\ZbbZbb{{
\unitlength=.07pt
\begin{picture}(192.00,120.00)(-16.00,0.00)
\bezier{20}(20.00,20.00)(30.00,10.00)(40.00,0.00)
\bezier{30}(80.00,80.00)(120.00,40.00)(160.00,0.00)
\bezier{30}(80.00,80.00)(40.00,40.00)(0.00,0.00)
\put(80.00,120.00){\line(0,-1){120.00}}
\end{picture}}
}

\def\bbZbbZ{{
\unitlength=0.07pt
\begin{picture}(192.00,120.00)(-16.00,0.00)
\bezier{20}(140.00,20.00)(130.00,10.00)(120.00,0.00)
\bezier{30}(80.00,80.00)(120.00,40.00)(160.00,0.00)
\bezier{30}(80.00,80.00)(40.00,40.00)(0.00,0.00)
\put(80.00,120.00){\line(0,-1){120.00}}
\end{picture}}
}

\def\bZbbZb{{
\unitlength=.07pt
\begin{picture}(192.00,120.00)(-16.00,0.00)
\bezier{20}(80.00,20.00)(90.00,10.00)(100.00,0.00)
\bezier{20}(80.00,20.00)(70.00,10.00)(60.00,0.00)
\put(80.00,40.00){\line(0,-1){20.00}}
\put(80.00,40.00){\line(0,1){0.00}}
\put(80.00,80.00){\line(0,-1){40.00}}
\put(80.00,120.00){\line(0,-1){40.00}}
\bezier{30}(80.00,80.00)(120.00,40.00)(160.00,0.00)
\bezier{30}(80.00,80.00)(40.00,40.00)(0.00,0.00)
\end{picture}}
}

\def\ZbbbZb{{
\unitlength=.07pt
\begin{picture}(192.00,120.00)(-16.00,0.00)
\put(40.00,0.00){\line(0,1){20.00}}
\put(40.00,40.00){\line(0,-1){40.00}}
\bezier{20}(40.00,40.00)(60.00,20.00)(80.00,0.00)
\put(80.00,40.00){\line(0,1){0.00}}
\put(80.00,120.00){\line(0,-1){40.00}}
\bezier{30}(80.00,80.00)(120.00,40.00)(160.00,0.00)
\bezier{30}(80.00,80.00)(40.00,40.00)(0.00,0.00)
\end{picture}}
}

\def\bZbbbZ{{
\unitlength=0.07pt
\begin{picture}(192.00,120.00)(-16.00,0.00)
\put(120.00,40.00){\line(0,-1){40.00}}
\bezier{20}(120.00,40.00)(100.00,20.00)(80.00,0.00)
\put(80.00,40.00){\line(0,1){0.00}}
\put(80.00,120.00){\line(0,-1){40.00}}
\bezier{30}(80.00,80.00)(120.00,40.00)(160.00,0.00)
\bezier{30}(80.00,80.00)(40.00,40.00)(0.00,0.00)
\end{picture}}
}

\def\ot{\otimes} 
\def\rada#1#2{#1,\ldots,#2}

\def\Rada#1#2#3{#1_{#2},\dots,#1_{#3}}
\def\landRada#1#2#3{#1_{#2} \land \dots \land #1_{#3}}
\def\pa{\partial}
\def\znamenko#1{{(-1)^{#1}}}
\def\doubless#1#2{{
\def\arraystretch{.5}
\begin{array}{c}
\mbox{\scriptsize $\scriptstyle #1$}
\\
\mbox{\scriptsize $\scriptstyle #2$}
\end{array}\def\arraystretch{1}
}}
\def\qed{\hspace*{\fill}
\mbox{\hphantom{mm}\rule{0.25cm}{0.25cm}}\\}

\def\btb{{\bullet \to \bullet}}     \def\max{{\overline b}}
\def\id{1\!\!1}                     \def\min{{\underline b}}
\def\bfk{{\bf k}} \def\maxx#1{#1_{\it max}} \def\minn#1{#1_{\it min}}
\def\ainfty{$A_\infty$}
\def\opinf{\underline{\cal A}}
\def\Ass{\mbox{\underline{${\cal A}${\it ss}}}}
\def\Asssoft{\underline{{\cal A}{\it ss}}}
\def\Deltasu{\Delta^{{\tt \hskip -1mm su}}\hskip -.1mm}

\def\Astrict{\hbox{${\tt strA}_\infty$}}
\def\Aweak{\hbox{${\tt shA}_\infty$}}

\def\leftrubber#1{\left(\rule{0mm}{#1mm}\right.}
\def\rightrubber#1{\left.\rule{0mm}{#1mm}\right)}
\def\br{{\bf R}}
\def\cc{{\cal C}}
\def\ck{{\cal K}}
\def\uck{\underline\ck}\def\normalK{K}\def\normalW{W}
\def\cw{{\cal W}}
\def\ucw{\underline\cw}

\def\pa{\partial}
\def\mbn{{\it mBin}(n)}
\def\into{
   {\hskip .7mm
   \unitlength=.14pt
   \begin{picture}(30.00,60.00)(0.00,0.00)
   \thicklines
   \put(30.00,0.00){\line(0,1){60.00}}
   \put(0.00,0.00){\line(1,0){30.00}}
   \end{picture}
   \hskip 1mm}
}

\title{Associahedra, cellular $W$-construction and
       products of $A_\infty$-algebras}

\author{Martin Markl\thanks{Supported by the grant GA \v CR 1019203.}
\hglue 4mm and
Steve Shnider\thanks{Supported by the Israel Academy of Sciences.}%
}
\begin{document}

\bibliographystyle{plain}

\maketitle

\begin{abstract}
  The aim of this paper is to construct a functorial tensor product of
  $A_\infty$-algebras or, equivalently, an explicit diagonal for the
  operad of cellular chains, over the integers, of the Stasheff
  associahedron. These construction were in fact already indicated
  in~\cite{umble-saneblidze}; we will try to give a more satisfactory
  presentation. We also prove that there does not exist
  a co-associative diagonal.
\end{abstract}

{\small
\vskip 3mm
\noindent
{\bf Table of content:} \ref{intro}.
                     Introduction
                     -- page~\pageref{intro}
\hfill\break\noindent
\hphantom{{\bf Table of content:\hskip .5mm}}  \ref{catprops}.
                     Categorial properties of diagonals and tensor products
                     -- page~\pageref{catprops}
\hfill\break\noindent
\hphantom{{\bf Table of content:\hskip .5mm}} \ref{dalsi_dva_roky}.
                     Calculus of oriented cell complexes of $K_n$ and $W_n$
                     -- page~\pageref{dalsi_dva_roky}
\hfill\break\noindent
\hphantom{{\bf Table of content:\hskip .5mm}}  \ref{maps}.
                     The chain maps $p$ and $q$
                     -- page~\pageref{maps}
\hfill\break\noindent
\hphantom{{\bf Table of content:\hskip .5mm}}  \ref{ElisAbeth}.
                     The Saneblidze-Umble diagonal
                     -- page~\pageref{ElisAbeth}
\hfill\break\noindent
\hphantom{{\bf Table of content:\hskip .5mm}}  \ref{tyden-pred-prohlidkou}.
                     Non-existence of a co-associative diagonal
                     -- page~\pageref{tyden-pred-prohlidkou}
\hfill\break\noindent
\hphantom{{\bf Table of content:\hskip .5mm}}  \ref{postponed}.
                     Remaining proofs
                     -- page~\pageref{postponed}
}

\section{Introduction}
\label{intro}

In this paper we study tensor products of
\ainfty-algebras. More precisely, given two
\ainfty-algebras $A = (V,\pa^V,\mu_2^V,\mu_3^V,\ldots)$ and $B=
(W,\pa^W,\mu_2^W,\mu_3^W,\ldots)$, we will be looking
for a functorial definition of an
\ainfty-structure $A\odot B$ that would extend the standard (non-associative)
dg-algebra structure on the tensor product $A \ot B$. This means that
the \ainfty-algebra $A \odot B$ will be of the form
$(V \otimes W, \pa,\mu_2,\mu_3,\ldots)$,
where $\pa$ is the usual differential on the tensor product,
\begin{equation}
\label{1}
\pa(v \ot w) := \pa^V(v) \ot w + \znamenko{\deg{v}} v \ot \pa^W(w)
\end{equation}
and the bilinear product $\mu_2$ is given by another standard formula
\begin{equation}
\mu_2(v' \ot w', v'' \ot w'') :=
\znamenko{\deg{v''}\deg{w'}}\mu^V(v',v'') \ot \mu^W(w',w''),
\end{equation}
where $v,v',v'' \in V$ and $w,w',w'' \in W$.

A ``coordinate-free'' formulation of the problem is the following. Let
$\opinf$ be the non-$\Sigma$ operad describing \ainfty-algebras
 (see~\cite[page~45]{markl-shnider-stasheff:book}),
that is, the minimal model of the non-$\Sigma$
operad $\Ass$ for associative algebras.
The above product is equivalent to a morphism of
dg-operads  (a {\em diagonal\/}) $\Delta : \opinf \to \opinf \ot \opinf$
such that $\Delta$ induces the usual diagonal $\Delta_{\Asssoft}$ on the
non-$\Sigma$ operad $\Ass = H_*(\opinf)$.

The existence of such a diagonal is not surprising
and follows from properties of minimal models for operads,
see~\cite[Proposition~3.136]{markl-shnider-stasheff:book}.
On the other hand, there is no way to
control the co-associativity of diagonals
constructed using this general argument and we will see below,
in Theorem~\ref{ElisaBeth},
that there, surprisingly enough,
{\em does not exist\/} a co-associative diagonal.

For practical purposes, such as applications in open string
theory~\cite{gaberdiel-zwiebach:NP97}, one needs
a tensor product (and therefore also a diagonal) given by an
explicit formula. Such an explicit diagonal was constructed by Umble
and Saneblidze in~\cite{umble-saneblidze}.
Our work was in fact motivated by our
unsuccessful attempts to understand their paper.
We will denote this diagonal by
$\Deltasu$ and call it the {\em  SU-diagonal\/}.
In this article we recall the definition of this diagonal and give a
conceptual explanation why it is well-defined.
The operad $\opinf$ can be identified with the operad of the
cellular  chain complexes of the
non-$\Sigma$ operad of associahedra, $\opinf\cong\cc_*(\uck)$
 (see~\cite[page~45]{markl-shnider-stasheff:book}), therefore, the required
diagonal is given by a family of chain maps

\[
\Delta_{\normalK_n}:
\cc_*(\normalK_n)\rightarrow\cc_*(\normalK_n)\ot\cc_*(\normalK_n),\ n
\geq 1,
\]
commuting with the induced operad structures and such
that $H_*(\Delta_{\normalK_n}) = \Delta_{\Asssoft}$.

The cells of the associahedra are not conducive to the definition of a
diagonal. There is, however, a cubical decomposition of the
associahedra provided by the {\em W-construction\/} of Boardman and
Vogt~\cite{boardman-vogt:73}, which is a homotopically equivalent
non-$\Sigma$ operad $\ucw = \{\normalW_n\}_{n \geq 1}$, for which
there is a canonical diagonal
\[
\Delta_{\normalW_n}: \cc_*(\normalW_n)\longrightarrow
\cc_*(\normalW_n)\ot\cc_*(\normalW_n),\ n\geq 1,
\]
induced by the cubical structure
(see~(\ref{pristi_utery_prohlidka})). A suitable diagonal on the
associahedra can be then obtained by transfering $\Delta_{\normalW_n}$
from $\ucw$ to $\uck$. More precisely, let
\begin{equation}
\label{trochu_se_stydim}
\cc_*(\normalW_n)\stackrel{p_n} \longrightarrow\cc_*(\normalK_n)
\stackrel{q_n}\longrightarrow \cc_*(\normalW_n),\ n \geq 1,
\end{equation}
be arbitrary operadic maps such that $H_*(p_n)$ and $H(q_n)$ are
identity endomorphisms of $\Ass(n)$, via the canonical identifications
\[
H_*(\cc_*(\normalW_n)) \cong \Ass(n) \cong H_*(\cc_*(\normalK_n)),\ n\geq 1.
\]
Then the formula
\begin{equation}
\label{zitra_prohlidka}
\Delta_{\normalK_n} :=(p_n\ot p_n)\circ\Delta_{\normalW_n}\circ q_n
\end{equation}
clearly defines a diagonal. In fact, it can be proved that the
operadic maps $p = \{p_n\}_{n \geq 1}$ and $q = \{q_n\}_{n \geq 1}$
with the above properties are homotopy inverses, but we will not need
this statement.

It remains to find maps in~(\ref{trochu_se_stydim}).
While there is an obvious and simple
definition of $q_n$, finding a suitable formula for $p_n$ is much
less obvious. We give an explicit and very natural definition inspired
by a formula in~\cite{umble-saneblidze}.

We will see that the operad of cellular chains
$\cc_*(\underline{W})$ can be described in terms of metric trees.
Similar {\em cellular $W$-constructions\/} on a given dg-operad
were considered by Kontsevich and Soibelman in~\cite{ko-so}. In
this terminology, the chain maps $p$ and $q$ are explicit homotopy
equivalences, defined over the integers, between the chain $W$-%
construction on the operad $\Ass$ and the minimal model $\opinf$
of $\Ass$, which give rise to explicit equivalences of the
categories of algebras over these dg-operads.

\section{Categorial properties of diagonals and tensor products}
\label{catprops}

Recall~\cite{markl:JPAA92} that there are two notions of
morphisms of \ainfty-algebras.
A {\em strict morphism\/} of \ainfty-algebras $(X,\pa,\mu_2,\mu_3,\ldots)$ and
$(Y,\pa,\nu_2,\nu_3,\ldots)$ is a linear map $f : X \to Y$
that commutes with all structure operations. A weaker notion
is that of a strongly homotopy (sh) morphism, given by a
sequence of maps $f_n : X^{\otimes n} \to Y$, $n \geq 1$, satisfying
rather complicated set of axioms (see, for
example,~\cite{markl:JPAA92,markl:ha}). Such a map is invertible if and only if
$f_1 : X \to Y$ is an isomorphism.
We will denote by \Astrict\ the category of \ainfty-algebras and their
strict morphisms, and \Aweak\ the category of
\ainfty-algebras and their sh{} morphisms.

As proved in~\cite[Proposition~3.136]{markl-shnider-stasheff:book},
any two diagonals $\Delta', \Delta'' : \opinf \to \opinf \ot \opinf$
are homotopic as maps of operads. Let $\odot'$ (resp.~$\odot''$)
denotes the tensor product induced by $\odot'$
(resp.~$\odot''$). Although $A \odot' B$ and $A \odot'' B$ are, in
general, not strictly isomorphic, the homotopy between  $\Delta'$ and
$\Delta''$ can be shown to induce a strongly homotopy
isomorphism between  $A \odot' B$ and $A \odot'' B$. Therefore we
obtain the following uniqueness:

\begin{proposition}
For any two \ainfty-algebras $A$, $B$, the \ainfty-algebras $A
\odot' B$ and $A \odot'' B$ are isomorphic in \Aweak.
\end{proposition}

We will prove in Theorem~\ref{ElisaBeth} that there are no co-associative
diagonals. This means that in general
\[
A\odot (B \odot C)  \not\cong (A \odot B) \odot C
\]
in the `strict' category \Astrict. On the other hand, as argued
in~\cite[Proposition~3.136]{markl-shnider-stasheff:book}, each diagonal
$\Delta$ is {\em homotopy associative\/} in the sense that the maps
$(\Delta \ot \id)\Delta$ and $(\id \ot \Delta)\Delta$ are
homotopic maps of operads, from which we infer:

\begin{proposition}
For any three \ainfty-algebras $A$, $B$ and $C$,
\[
A\odot (B \odot C)  \cong (A \odot B) \odot C
\]
in the `weak' category \Aweak.
\end{proposition}

By the same argument, one can also prove

\begin{proposition}
For any two \ainfty-algebras $A$ and $B$,
\[
A\odot B \cong B \odot A
\]
in \Aweak.
\end{proposition}

This naturally rises the question whether \Aweak\ with a product
$\odot$ based on an appropriate diagonal is a (possibly
symmetric) monoidal category.
Even to formulate this question precisely, one more step should be
completed.

While it is clear that $\odot$ is a functor $\Astrict \times
\Astrict \to \Astrict$, to make it a functor $\Aweak \times
\Aweak \to \Aweak$, one should define, for two sh{} morphisms $f: A'
\to A''$ and $g : B' \to B''$, a `product' $f \odot g : A' \odot
A'' \to B' \odot B''$. One should then consider a functorial
`associator' $\Phi_{A,B,C} : A\odot (B \odot C)
\to  (A \odot B) \odot C$ and a `symmetry'
$\sigma_{A,B} : A\odot B \to B \odot A$.

The above objects exist by general nonsense, but it is not clear
whether they fulfill the axioms of a (symmetric) monoidal category
(the pentagon and the hexagons), although it is quite possible that
for some special choices of the above data these axioms are satisfied.
On a more abstract level, the `full' functorial monoidal product $A,B
\mapsto A \odot B$, $f,g, \mapsto f \odot g$ in $\Aweak$ means to
construct a `diagonal' in the minimal model of the two-colored operad
$\Ass_\btb$ describing homomorphisms of associative algebras,
satisfying some additional properties which do not follow from a
general nonsense.

\section{Calculus of oriented cell complexes of $K_n$ and $W_n$}
\label{dalsi_dva_roky}

All operads ${\mathcal P}$ considered in this paper are such
that ${\mathcal P}(0)$ is trivial and that ${\mathcal P}(1)$ is
isomorphic to the ground field. The category of operads with this
property is equivalent to the category of {\em pseudo-operads\/}
${\mathcal P}$ such that ${\mathcal P}(0) = {\mathcal P}(1) =
0$, the equivalence being given by forgetting the $n=1$ piece.
This, roughly speaking, means that we may ignore operadic units,
see~\cite[Observation~1.2]{markl:zebrulka} for
details. {\em Therefore, for the rest of this paper, an operad means a
pseudo-operad with ${\mathcal P}(0) = {\mathcal P}(1) = 0$.}

First, we establish some notation. Let $\uck=\{\normalK_{n}\}_{n\geq
2}$ be the non-$\Sigma$ operad of associahedra. The topological cell
complex $\normalK_n$ can be realized as a convex polytope in
$\br^{n-2}$, with $k$-cells labeled by the planar rooted trees with
$n$ leaves and $n-k-2$ internal edges, or equivalently by
$(n-k-2)$-fold bracketings of $n$ elements,
see~\cite[II.1.6]{markl-shnider-stasheff:book}. For example, $0$-cells
correspond to binary trees with $n$ leaves, or equivalently, full
bracketings of $n$ elements.  All our constructions will be expressed
in terms of rooted planar trees although there is clearly an
underlying geometric meaning based on the polytope realization of
$\normalK_n$. Boardman and Vogt have defined
in~\cite{boardman-vogt:73} a cubical subdivision of the cells of
$\normalK_n$, for $n\geq 2$, giving rise to a cubical cell complex
known as the $W$-construction, $\normalW_n$.  See Figure~6
of~\cite[Section~II.2.8]{markl-shnider-stasheff:book} for $\normalW_4$
represented as a cubical subdivision of $\normalK_4$.

The cells of $\normalW_n$ are in one-to-one
correspondence with ``metric $n$-trees,'' that is, planar rooted
trees with $n$ leaves and with internal edges labeled either
``metric'' or ``non-metric.'' The metric $n$-trees with $k$
metric edges label the topological $k$-cells of $\normalW_n$.
A cubical cell is called  {\em an interior cell\/} if the
labeling tree has only metric edges.  In the geometric realization
the interior cells are in the interior of the convex polytope.

In order to define the boundary operators on the  complexes $\cc_*(\uck):
=\{\cc_*(\normalK_n)\}_{n\geq 2}$ and
$\cc_*(\ucw):=\{\cc_*(\normalW_n)\}_{n\geq 2}$
(non-$\Sigma$ operads in the category of chain complexes),
we have to introduce an orientation on the cells.
Let $T$ be a planar rooted tree with internal edges labeled $e_1,\ldots, e_m$.
Two orderings  $e_{i_1},\ldots,e_{i_m}$  and $e_{j_1},\ldots,e_{j_m}$
will be called {\em equivalent\/} if
they are related  by an even permutation.
The equivalence class corresponding to an ordering  $e_{i_1},\ldots,e_{i_m}$
will be called an orientation and denoted
 $e_{i_1}\wedge \cdots \wedge e_{i_m}$.

\begin{definition}
An oriented $k$-cell in $\normalK_n$ is a pair  $(T, \omega)$ where
$T$ is a planar rooted tree with $n$ leaves and $n-k-2$ internal edges
and $\omega$ is an orientation.
Let $\cc_k(\normalK_n)$ be the vector space spanned by the oriented $k$-cells
in $\normalK_n$ modulo the relation $(T,\omega)=-(T,\omega')$ where
$\omega$ and
$\omega'$ are the two distinct orientations.

An oriented metric $k$-cell in $\normalW_n$ is a pair  $(T, \omega)$ where
$T$ is a metric tree with $n$ leaves and $k$  metric edges and
$\omega$ is an orientation of the metric edges. Let $\cc_k(\normalW_n)$
be the vector space spanned by the oriented $k$-cells in $\normalW_n$ modulo
 the relation $(T,\omega)=-(T, \omega')$
where $\omega$ and $\omega'$ are the two distinct orientations.
\end{definition}

The operad composition law
$$\circ_i:\cc_k(\normalK_r)\ot \cc_l(\normalK_s)\longrightarrow
\cc_{k+l}(\normalK_{r+s-1})$$
is  defined on the basis elements by
\begin{equation}\label{circick}
(T,\omega)\circ_i (T',\omega'):=(-1)^{rl +i(s+1)}(T\circ_i T',
\omega\wedge \omega'\wedge e),
\end{equation}
where $\circ_i$ is defined on planar rooted trees in the standard way,
grafting the second tree onto the $i$-th leaf of the first, and
$\omega\wedge \omega'\wedge e$ is the
concatenation of the two orientations, with the new edge created by
grafting labeled $e$.
The operad composition law
$$
\circ_i:\cc_k(\normalW_r)\ot\cc_l(\normalW_s)\longrightarrow
\cc_{k+l}(\normalW_{r+s-1})
$$
is  defined on the basis elements by
\begin{equation}
\label{circicw}
(T,\omega)\circ_i (T',\omega'):=(T\circ_i T', \omega\wedge\omega'),
\end{equation}

A heuristic  explanation of why we don't need
any signs in the above display is that the orientation of the cells of
$\ucw$ defined in terms of metric edges is geometric in the sense that
the number of metric edges is the same as the dimension of the cell.
In the case of $\cc_*(\ucw)$ the new edge created by the
grafting is non-metric and so does not appear in the ordering of
metric vertices.

The boundary operator on $\cc_*(\normalK_n)$ is defined by
\begin{equation}
\pa_\normalK (T,e_1\wedge\cdots\wedge e_m):=
\sum_{\{T'\, | \,T'/{e'}=T\}} (T',e'\wedge e_1\wedge\cdots\wedge e_m),
\end{equation}
where the sum is over all trees $T'$ with an edge $e'$,
such that when $e'$ is collapsed $T'$ reduces to $T$.
The condition $\pa_\normalK^2=0$ follows immediately from the identities
$$
(T'',\, e'\wedge e''\wedge e_1\wedge \cdots\wedge e_m)
=-(T'',\,e''\wedge e'\wedge e_1\wedge\cdots\wedge e_m).
$$

Next we define the boundary operator on the complex  $\cc_*(\normalW_n)$.
Let $T$ be a metric tree and  $ e_1\wedge\cdots\wedge e_k$ an
orientation:
\begin{eqnarray}
\label{pacw}
\lefteqn{
\pa_\normalW(T, e_1\wedge\cdots\wedge e_k):=}
\\
\nonumber
&&
\sum_{1 \leq i \leq k}(-1)^{i}
\left[
(T/e_i, e_1\wedge\cdots \hat{e_i} \cdots\wedge e_k) -
(T_i,\, e_1\wedge \cdots\hat{e_i} \cdots\wedge e_k)
\right],
\end{eqnarray}
where $T_i$ is the same (unlabeled) tree as $T$ but with the metric edge
$e_i$  changed to a non-metric edge. As above, the condition
$\pa_\normalW^2=0$ follows from the relations for the orientation elements.

In the rest of this section we introduce `standard orientations' for top
dimensional cells of $W_n$ and $0$-dimensional cells of $K_n$.
There is a partial order relation on rooted planar binary trees
given by the associator which moves a vertex
to the right and changes the outgoing edge from a right leaning position to
a left-leaning position, as shown in Figure~\ref{arrows}.
\begin{figure}[t]
\begin{center}
{
\unitlength=1pt
\begin{picture}(210.00,50.00)(-30.00,10.00)
\thicklines
%
\put(160.00,30.00){\vector(-1,0){50.00}}
%
\put(20,0){
\put(60.00,60.00){\makebox(0.00,0.00){$\bullet((\bullet\bullet)\bullet)$}}
\put(60.00,10.00){\makebox(0.00,0.00){$\bullet$}}
\put(50.00,10.00){\makebox(0.00,0.00){$\beta$}}
\put(70.00,20.00){\makebox(0.00,0.00){$\bullet$}}
\put(80.00,20.00){\makebox(0.00,0.00){$\alpha$}}
\put(60.00,30.00){\makebox(0.00,0.00){$\bullet$}}
\put(60.00,10.00){\line(1,-1){10.00}}
\put(70.00,20.00){\line(-1,-1){20.00}}
\put(60.00,30.00){\line(1,-1){30.00}}
\put(60.00,30.00){\line(-1,-1){30.00}}
\put(60.00,50.00){\line(0,-1){20.00}}
}
%
\put(50.00,30.00){\vector(-1,0){40.00}}
%
\put(-30.00,60.00){\makebox(0.00,0.00){$\bullet(\bullet(\bullet\bullet))$}}
\put(-10.00,10.00){\makebox(0.00,0.00){$\bullet$}}
\put(-20.00,20.00){\makebox(0.00,0.00){$\bullet$}}
\put(0.00,10.00){\makebox(0.00,0.00){$\beta$}}
\put(-30.00,30.00){\makebox(0.00,0.00){$\bullet$}}
\put(-10.00,10.00){\line(-1,-1){10.00}}
\put(-20.00,20.00){\line(-1,-1){20.00}}
\put(-30.00,30.00){\line(1,-1){30.00}}
\put(-30.00,30.00){\line(-1,-1){30.00}}
\put(-30.00,50.00){\line(0,-1){20.00}}
%
%
\put(10,0){
\put(180.00,60.00){\makebox(0.00,0.00){$(\bullet(\bullet\bullet))\bullet$}}
\put(180.00,10.00){\makebox(0.00,0.00){$\bullet$}}
\put(170.00,20.00){\makebox(0.00,0.00){$\bullet$}}
\put(160.00,20.00){\makebox(0.00,0.00){$\alpha$}}
\put(180.00,30.00){\makebox(0.00,0.00){$\bullet$}}
\put(180.00,10.00){\line(-1,-1){10.00}}
\put(170.00,20.00){\line(1,-1){20.00}}
\put(180.00,30.00){\line(1,-1){30.00}}
\put(180.00,30.00){\line(-1,-1){30.00}}
\put(180.00,50.00){\line(0,-1){20.00}}
}
\put(25.00,60.00){\makebox(0.00,0.00){\large $>$}}
\put(135.00,60.00){\makebox(0.00,0.00){\large $>$}}
\end{picture}}
\end{center}
\caption{%
\label{arrows}%
The partial order on the set of binary trees. The first on the right
arrow moves the vertex $\alpha$ and the second arrow moves the
vertex $\beta$. }
\end{figure}
The {\em standard orientation $\omega_{\max(n)}$} of the maximal
fully metric binary tree $\max(n)$ (all internal edges leaning to
the left) is given by enumerating the internal edges in sequence,
starting with $e_1$, the edge adjacent to the root, and continuing
$e_2,\ldots, e_{n-2}$ in sequence going away from the root, see
Figure~\ref{houpacka}.

\begin{figure}
\centering
\thicklines
\unitlength 4mm
\begin{picture}(16,9)(2,2)
\put(9,9){\line(1,-1){7}}
\put(9,9){\line(-1,-1){7}}
\put(10,8){\line(-1,-1){6}}
\put(11,7){\line(-1,-1){5}}
\put(13,5){\line(-1,-1){3}}
\put(14,4){\line(-1,-1){2}}
\put(15,3){\line(-1,-1){1}}
\put(-.5,-.5){
\put(10,9){\makebox(0,0)[lb]{$e_1$}}
\put(11,8){\makebox(0,0)[lb]{$e_2$}}
\put(12.5,6.5){\makebox(0,0)[lb]{\scriptsize $\ddots$}}
\put(14,5){\makebox(0,0)[lb]{$e_{n-3}$}}
\put(15,4){\makebox(0,0)[lb]{$e_{n-2}$}}}
\end{picture}
\caption{The maximal binary tree $\max(n)$.}
\label{houpacka}
\end{figure}
\begin{figure}[t]
\begin{center}
\setlength{\unitlength}{2500sp}%
\begin{picture}(6000,5649)(0,-4873)
\thicklines
%
\put(2050,-4200){\vector(1,0){3050}}
\put(1300,-1900){\vector(1,-3){760}}
\put(5100,-4200){\vector(1,3){760}}
\put(1300,-1900){\vector(4,3){2280}}
\put(3600,-200){\vector( 4,-3){2280}}
%
\put(-1800,-2000){
\makebox(0,0){\put(-1300,0){$\left(\rule{0mm}{5mm}\right.$}
\setlength{\unitlength}{.2cm}
\thinlines
\put(-8,-9){
\bezier{200}(1,8)(2.5,9.5)(4,11)
\bezier{200}(4,11)(5.5,9.5)(7,8)
\bezier{100}(3,10)(4,9)(5,8)
\bezier{50}(2,9)(2.5,8.5)(3,8)
\put(3.5,10.5){\makebox(0,0)[rb]{\scriptsize $e_1$}}
\put(2.5,9.5){\makebox(0,0)[rb]{\scriptsize $e_2$}}
}
$,-e_1 \land e_2
\left.\rule{0mm}{5mm}\right) =$
}
\put(1000,-700){
\makebox(0,0){\put(-1300,0){$\left(\rule{0mm}{5mm}\right.$}
\setlength{\unitlength}{.2cm}
\thinlines
\put(-8,-9){
\bezier{200}(1,8)(2.5,9.5)(4,11)
\bezier{200}(4,11)(5.5,9.5)(7,8)
\bezier{100}(3,10)(4,9)(5,8)
\bezier{50}(2,9)(2.5,8.5)(3,8)
\put(3.5,10.5){\makebox(0,0)[rb]{\scriptsize $e_2$}}
\put(2.5,9.5){\makebox(0,0)[rb]{\scriptsize $e_1$}}
}
$,+e_1 \land e_2
\left.\rule{0mm}{5mm}\right)$
}}}
\put(1300,-1900){\makebox(0,0){$\bullet$}}

\put(1800,-5200){\makebox(0,0){\put(-1300,0){$\left(\rule{0mm}{5mm}\right.$}
\setlength{\unitlength}{.2cm}
\thinlines
\put(-15,-9){
\bezier{200}(8,8)(9.5,9.5)(11,11)
\bezier{200}(11,11)(12.5,9.5)(14,8)
\bezier{100}(9.5,9.5)(10.5,8.5)(11,8)
\bezier{100}(10.5,8.5)(10.25,8.25)(10,8)
\put(10.5,10.5){\makebox(0,0)[rb]{\scriptsize $e_1$}}
\put(10,9){\makebox(0,0)[lb]{\scriptsize $e_2$}}
}
$,+e_1 \land e_2
\left.\rule{0mm}{5mm}\right)$
}}
\put(2060,-4200){\makebox(0,0){$\bullet$}}

\put(5100,-5200){\makebox(0,0){\put(-1300,0){$\left(\rule{0mm}{5mm}\right.$}
\setlength{\unitlength}{.2cm}
\thinlines
\put(-5,-9){
\bezier{200}(4,8)(2.5,9.5)(1,11)
\bezier{200}(1,11)(-0.5,9.5)(-2,8)
\bezier{100}(2.5,9.5)(1.5,8.5)(1,8)
\bezier{100}(1.5,8.5)(1.75,8.25)(2,8)
\put(1.5,10.5){\makebox(0,0)[lb]{\scriptsize $e_1$}}
\put(2.1,9.1){\makebox(0,0)[rb]{\scriptsize $e_2$}}
}
$,-e_1 \land e_2
\left.\rule{0mm}{5mm}\right)$
}}
\put(5100,-4200){\makebox(0,0){$\bullet$}}
%
\put(7500,-2400){\makebox(0,0){\put(-1300,0){$\left(\rule{0mm}{5mm}\right.$}
\setlength{\unitlength}{.2cm}
\thinlines
\put(-5,-9){
\bezier{200}(4,8)(2.5,9.5)(1,11)
\bezier{200}(1,11)(-0.5,9.5)(-2,8)
\bezier{100}(2,10)(1,9)(0,8)
\bezier{50}(3,9)(2.5,8.5)(2,8)
\put(1.5,10.5){\makebox(0,0)[lb]{\scriptsize $e_1$}}
\put(2.5,9.5){\makebox(0,0)[lb]{\scriptsize $e_2$}}
}
$,+e_1 \land e_2
\left.\rule{0mm}{5mm}\right)$
}}
\put(5850,-1900){\makebox(0,0){$\bullet$}}
%
\put(3500,-200){\makebox(0,0){\put(-1300,0){$\left(\rule{0mm}{5mm}\right.$}
\setlength{\unitlength}{.2cm}
\thinlines
\put(-5,-9){
\bezier{200}(4,8)(2.5,9.5)(1,11)
\bezier{200}(1,11)(-0.5,9.5)(-2,8)
\bezier{50}(-1,9)(-.5,8.5)(0,8)
\bezier{50}(3,9)(2.5,8.5)(2,8)
\put(0,10){\makebox(0,0)[rb]{\scriptsize $e_1$}}
\put(2.5,10){\makebox(0,0)[lb]{\scriptsize $e_2$}}
}
$,-e_1 \land e_2
\left.\rule{0mm}{5mm}\right)$
}}
\put(3600,-200){\makebox(0,0){$\bullet$}}
%
\end{picture}%
\end{center}
\caption { The rule for defining the standard orientation of fully
 metric trees is
 illustrated for the pentagon in the figure above. This example also
 verifies that the definition is independent of the path. }
\label{stanori}
\end{figure}
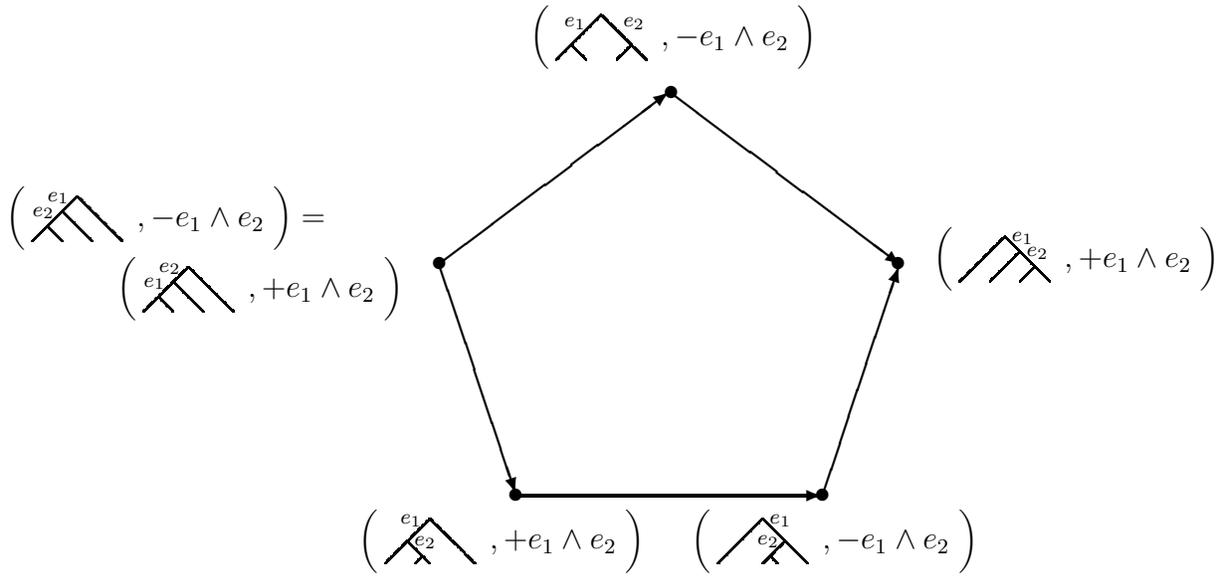

%
The standard orientation $\omega_T$ of a non-maximal fully metric
binary tree $T$ is determined by a sequence of sign changes and
relabelings along a path from $\max(n)$ to $T$ in the associahedron. See
Figure~\ref{stanori} for the standard orientations of binary trees
with four leaves.
To check that this rule gives and unambiguous definition of the orientation,
it is sufficient  (thanks to Mac~Lane's Coherence Theorem) to verify
that the definition is independent of path in  the pentagon (expressing
coherence of the associator) and in the square (expressing naturality).
The verification for the pentagon is given in Figure~\ref{stanori}.
The verification for the square is straightforward and
follows from the functoriality.
Therefore each fully metric binary $n$-tree $T$ together with its
standard orientation $\omega_T$ determines an element $(T,\omega_T)
\in \cc_{n-2}(\normalW_n)$.

We also define the {\em standard orientation\/} $\xi_T$ of a
binary $n$-tree $T$ representing a $0$-cell of $\cc_0(K_n)$
inductively as follows. The only binary $2$-tree representing a
$0$-cell of $\cc_0(K_2)$ has no internal edges, and its canonical
orientation is given by assigning the $+1$-sign to this tree. The
canonical orientation of any binary tree would be then determined
by the formula
\[
(S,\xi_S) \circ (T,\xi_T) = (S \circ_i T, \xi_{S \circ_i T}),
\]
once we checked that there was no ambiguity. This can be done
exactly as in the previous paragraph for $\omega_T$ . For example,
we immediately get the following standard orientations:
\[
\left(\unitlength 1mm
\linethickness{0.4pt}
\begin{picture}(10,4)(7,7)
\put(7,7){\line(1,1){4}}
\put(11,11){\line(1,-1){4}}
\put(15,7){\line(0,1){0}}
\qbezier(9,9)(10,8)(11,7)
\put(8.5,10.5){\makebox(0,0)[lb]{\mbox {\scriptsize $a$}}}
\end{picture}, -a
\right) \mbox{ and }
 \left( \unitlength 1mm \linethickness{0.4pt}
\begin{picture}(10,4)(-15.5,7)
\put(-7,7){\line(-1,1){4}}
\put(-11,11){\line(-1,-1){4}}
\put(-15,7){\line(0,1){0}}
\qbezier(-9,9)(-10,8)(-11,7)
\put(-8.5,10.5){\makebox(0,0){\mbox {\scriptsize $e$}}}
\end{picture}, e
\right).
\]
We recommend as an exercise to verify that the standard
orientation $\xi_{\max(n)}$ of the maximal binary tree in
Figure~\ref{stanori}, this time considered as a $0$ cell of $K_n$
is, for $n>2$,
\[
(\max(n),\xi_{\max(n)}) := (-1)^{(n-2)(n-3)/2} \cdot (\max(n),e_1
\land \cdots \land e_{n-2}) \in \cc_0(K_n)
\]
and that the standard orientation of the minimal binary $n$-tree
$\min(n)$ with the interior edges (all are right-leaning) enumerated
in sequence going away from the root, is given as
\[
(\min(n),\xi_{\min(n)}) := (-1)^{n} \cdot (\min(n),e_1
\land \cdots \land e_{n-2}) \in \cc_0(K_n).
\]

\section{The chain maps $p$ and $q$}
\label{maps}

The goal of this section is to construct maps $p : \cc_*(\uck) \to
\cc_*(\ucw)$ (Definition~\ref{propdefp}) and $q : \cc_*(\ucw) \to
\cc_*(\uck)$ (Definition~\ref{defq}) with the properties discussed in
Section~\ref{intro}. The proofs that that these maps are indeed chain
maps (Proposition~\ref{pisu_v_susarne} and
Proposition~\ref{ale_porad_na_ni_myslim}) are postponed to
Section~\ref{postponed}.

As an operad in the category of vector spaces, $\cc_*(\uck)$ is a free
operad generated by the collection with arity $n$ component, a
one-dimensional subspace concentrated in degree $n-2$ spanned by
corolla with $n$ leaves, and $\cc_*(\ucw)$ is a free operad generated
by the collection with arity $n$ component, the vector space with
basis the set of purely metric planar rooted trees with $n$ leaves.
Since a operadic map of a free operad is determined by its value on
generators, the operadic chain map
$\cc_*(\uck)\stackrel{q}\longrightarrow\cc_*(\ucw)$ is determined by
its value on corollae, and the operadic chain map
$\cc_*(\ucw)\stackrel{p}\longrightarrow\cc_*(\uck)$ is determined by
its value on purely metric trees.

 Let $c(n)$ be the corolla with $n$ leaves; since there are no
internal edges, we denote the orientation by the symbol $1$, and
adopt the convention that
$$1\wedge e_1\wedge\cdots\wedge e_k:=e_1\wedge\cdots\wedge e_k.$$

\begin{definition}
\label{defq}
Let ${\it mBin}(n)$ be the set of $n-2$ cells of $\normalW_n$
corresponding to the fully metric planar rooted binary trees with $n$
leaves and standard orientation.  Then $q(c(n),1)$ is defined as a
sum over ${\it mBin}(n)$:

\begin{equation}\label{qdef}
q(c(n),1):=\sum_{T\in\mbn}  (T,\omega_T).
\end{equation}
The operadic extension of $q$ to the free operad
$\cc(\uck)$ , which map will also be denoted $q$, defines a
morphism of operads in the category of graded vector spaces.
\end{definition}

\begin{proposition}
\label{ale_porad_na_ni_myslim}
The morphism $q$ described in Definition~\ref{defq}
commutes with the boundary operators,
\begin{equation}\label{qpa}
q(\pa_\normalK(c(n),1)=\pa_\normalW q(c(n),1),
\end{equation}
and therefore is a morphism of operads in the category of chain complexes.
\end{proposition}

The proof of Proposition~\ref{ale_porad_na_ni_myslim} is postponed
to Section~\ref{postponed}.  The operad chain map
$p:\cc_*(\ucw)\rightarrow \cc_*(\uck)$ is determined by its value
on fully metric trees. Before giving the precise definition, we
will give a conceptual description. As a topological cell complex,
the associahedron can be realized as a convex polytope
$\normalK_n\subset{\br^{n-2}}$.  The cubical cell complex
$\normalW_n$ is a decomposition of the associahedral $k$-cells
into $k$-cubes. The interior $k$-cell of $\normalW_n$ labeled by a
purely metric tree $T$ with $k$ edges is transverse to the $n-2-k$
cell of $\normalK_n$ labeled by the same tree.  Let $\minn{T}$ be
the binary tree labeling the minimal vertex of this transverse
cell in $\normalK_n$.

The image
$p(T)$ is defined as the sum with appropriate signs of all the
$k$-cells in $\cc_k(\normalK_n)$ all of whose vertices are labeled by
binary trees less than or equal to $\minn{T}$ relative to the partial
order on binary trees.

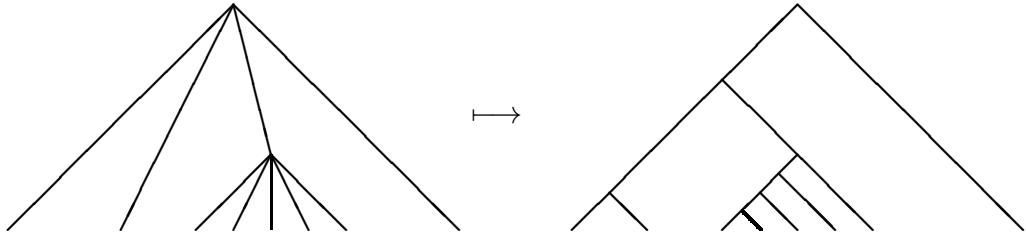
\begin{figure}
\begin{center}
\unitlength 5mm
\thicklines
\begin{picture}(29,6)(1,3)
\put(2,3){\line(1,1){6}}
\put(8,9){\line(1,-1){6}}
\put(8,9){\line(-1,-2){3}}
\put(9,5){\line(-1,-1){2}}
\put(9,5){\line(1,-1){2}}
\put(9,5){\line(1,-2){1}}
\put(9,5){\line(-1,-2){1}}
\put(9,5){\line(0,-1){2}}
\put(8,9){\line(1,-4){1}}
\put(17,3){\line(1,1){6}}
\put(23,9){\line(1,-1){6}}
\put(15,6){\makebox(0,0)[cc]{$\longmapsto$}}
\put(18,4){\line(1,-1){1}}
\put(21,7){\line(1,-1){4}}
\put(23,5){\line(-1,-1){2}}
\put(22,4){\line(1,-1){1}}
\put(22.5,4.5){\line(1,-1){1.5}}
\put(21.5,3.5){\bezier{50}(0,0)(.25,-.25)(.5,-.5)}

\end{picture}
\end{center}
\caption{An example of the filling-in procedure passing from a
fully metric tree
$T$ to the binary  tree $\minn{T}$.}
\label{bude_zitra_lat?}
\end{figure}

The tree $\minn{T}$ is created by ``filling-in'' the non-binary
vertices of $T$. A vertex in $T$ with $r$ input edges, $r>2$, is
replaced in $\minn{T}$ by the minimal binary tree with $r$ leaves,
which introduces $r-2$ new right-leaning edges.  When this procedure
is carried out at all the non-binary vertices of $T$, it adds $n-2-k$
new edges, all of which are right-leaning. See
Figure~\ref{bude_zitra_lat?} for an example of this procedure. In
exactly the same way, one defines $\maxx T$ as the binary tree obtained from
$T$ by filling-in the non-binary vertices by left-leaning edges.

In order for a binary tree $S$ to be the maximal vertex of a $k$-cell
in $\normalK_n$, it must contain at least $k$ left-leaning edges,
since an associativity move applied to a binary tree replaces a
right-leaning edge with a left-leaning edge (see Figure~\ref{arrows})
and the tree labeling the maximal vertex of $k$-cell is the output of
at least $k$ distinct associativity moves, corresponding to the $k$
one-cells of the associahedron which meet at the given vertex.

If $T$ is an interior $k$-cell in $\normalW_n$, then $\minn{T}$
cannot have more than $k$ left-leaning edges, since the new edges in
$\minn{T}$ are all right-leaning.  Since the number of left leaning
edges in a binary tree is a non-decreasing function relative to the
partial order, if $\minn{T}$ has less than $k$ left-leaning edges,
there are no $k$-cells less than $\minn{T}$ and we put $p(T) = 0$.

Given $(T,e_1\wedge\cdots\wedge e_k)\in\cc_k(\normalW_n)$, such that
$\minn{T}$ has $k$ left-leaning edges, then each edge $e_i$
corresponds to an edge in $\minn{T}$ which we also denote $e_i$.
Choose any labeling $f_1 \land \cdots \land f_{n-k-2}$ of the new
edges, and let $\xi_{\minn{T}}$ be the standard orientation of
$\minn{T}$ considered as the label for a $0$-cell of $K_n$,
$$
\xi_{\minn{T}}= \eta \cdot
e_1\wedge\cdots\wedge e_k\wedge f_1\wedge\cdots\wedge f_{n-k-2},\
\eta \in \{-1,+1\},
$$
then leading term of $p(T,e_1\wedge\cdots\wedge e_k)$ will be
\begin{eqnarray}
\label{Jsem_nemocny_nebo_je_mi_spatne_po_prasku?}
\lefteqn{
(\minn{T}/\{e_1,\ldots,e_k\},e_1\wedge\cdots\wedge e_k
\into \xi_{\minn{T}}) = \hskip 2cm}
\\
\nonumber
&&\hskip 2cm=
\eta \cdot (-1)^{k(k-1)/2}
(\minn{T}/\{e_1,\ldots,e_k\}, f_1 \land \cdots \land f_{n-k-2}).
\end{eqnarray}
In the above display, $\into$ is the
contraction  relative to the pairing $\langle e_i,e_j \rangle
:=\delta^i_j$.

For any binary tree $S<\minn{T}$ with $k$ left-leaning edges,
we will describe a method (analogous to the definition of the
standard orientation) of assigning in a  unique way
a labeling of the left-leaning edges by the labels $e_1,\ldots,e_k$.
First, we describe a rule which determines the  labeling of the
left-leaning edges in a tree given  the labeling of the left-leaning
edges in an adjacent tree (related by one associativity).
Consider a binary tree with labels only on the left-leaning edges,
adjacent trees
are related by replacing configuration of two edges
\unitlength .2cm
\begin{picture}(1.8,2)(-.3,.3)
\thicklines
\put(0,0){\bezier{30}(0,0)(.5,.5)(1,1)}
\put(1,1){\bezier{30}(0,0)(-.5,.5)(-1,1)}
\end{picture}
by the configuration
\unitlength .2cm
\begin{picture}(1.8,2)(-.8,.3)
\thicklines
\put(0,0){\bezier{30}(0,0)(-.5,.5)(-1,1)}
\put(-1,1){\bezier{30}(0,0)(.5,.5)(1,1)}
\end{picture}
(going from the greater tree  to the lesser tree). If both edges are
internal, the rule is simply to use the same label for the
left-leaning edge in both configurations. The ambiguity of the path
connecting two trees resolves into a sequence of pentagons and
squares and the validity of the
definition is checked by considering these two figures. On the other hand,
if the lower edge is a leaf, the new configuration has a right-leaning
edge in place of a left-leaning edge and the resulting binary tree
has less than $k$ left-leaning edges so that
there are no $k$-cells less than it and the contribution to $p(T)$ is zero.
For example, in Figure~\ref{arrows}, one of the
associativity moves preserves the number
of left-leaning edges and the other changes the number by one.

We can now give the full definition of $p$:

\begin{definition}
\label{propdefp}
Define a function on  oriented fully
metric trees $(T,e_1\wedge\cdots\wedge e_k)$ by
\begin{equation}\label{defp}
p(T,e_1\wedge\cdots\wedge e_k):=\sum (S/\{e_1,\ldots,e_k\},e_1\wedge
\cdots\wedge e_k\into\xi),
\end{equation}
where the sum is over binary trees $S$ less than or equal to
$\minn{T}$ with $k$ left-leaning edges labeled $e_1,\ldots,e_k$
according to the procedure described above, $\xi$ is the standard
orientation of the binary tree $S$ and $\into$ is the same
contraction as
in~(\ref{Jsem_nemocny_nebo_je_mi_spatne_po_prasku?}).

The function $p$ has a unique extension to a morphism (denoted also by
the same symbol)  $p:\cc_*(\normalW_*)\rightarrow \cc_*(\normalK_*)$
of operads in the category of graded vector spaces.
\end{definition}

\begin{exercise}
\label{osculetur_me}
{\rm
Verify that
\[
p(\max(n),\omega_{\max(n)})=(c(n),1)
\mbox { and }
p(c(n),1) = (\min(n),\xi_\min).
\]
Note that the first equation involves $(n-2)$-cells and the second
involves $0$-cells. Observe also that, modulo
orientations,~(\ref{defp}) is the sum of all trees $U$ with $n-k$
interior edges such that $\maxx U \leq \minn T$.
}
\end{exercise}
\begin{example}
\label{minuly_tyden_jsem_byl_s_Jitkou}
{\rm\
Let us describe explicitly the map $p : \cc_*(\normalW_n) \to
\cc_*(\normalK_n)$ for some small $n$. For $n=1$ and $2$, $p$ is given by
\[
p\left(
\unitlength1pt
\begin{picture}(6,12)(-1.5,1)
\put(0,0){\line(0,1){10}}
\end{picture},1
\right)
:=
\left(
\unitlength1pt
\begin{picture}(6,12)(-1.5,1)
\put(0,0){\line(0,1){10}}
\end{picture},1
\right)
\mbox { and }
p\left(
\unitlength.25cm
\begin{picture}(2,1)
\put(0,0){\bezier{30}(0,0)(.5,.5)(1,1)}
\put(1,0){\bezier{30}(0,1)(.5,.5)(1,0)}
\end{picture}\ , 1
\right) :=
\left(
\unitlength.25cm
\begin{picture}(2,1)
\put(0,0){\bezier{30}(0,0)(.5,.5)(1,1)}
\put(1,0){\bezier{30}(0,1)(.5,.5)(1,0)}
\end{picture}\ , 1
\right).
\]
For $n=3$,
\[
p\left(
\unitlength.4cm
\begin{picture}(2,1)
\put(0,0){\bezier{50}(0,0)(.5,.5)(1,1)}
\put(1,0){\bezier{50}(0,1)(.5,.5)(1,0)}
\put(1,0){\line(0,1){1}}
\end{picture}\ , 1
\right) :=
\left(
\unitlength 1mm
\linethickness{0.4pt}
\begin{picture}(10,4)(7,7)
\put(7,7){\line(1,1){4}}
\put(11,11){\line(1,-1){4}}
\put(15,7){\line(0,1){0}}
\qbezier(9,9)(10,8)(11,7)
\put(8.5,10.5){\makebox(0,0)[lb]{\mbox {\scriptsize $a$}}}
\end{picture}, -a
\right),\
p\left(
\unitlength 1mm
\linethickness{0.4pt}
\begin{picture}(10,4)(-15.5,7)
\put(-7,7){\line(-1,1){4}}
\put(-11,11){\line(-1,-1){4}}
\put(-15,7){\line(0,1){0}}
\qbezier(-9,9)(-10,8)(-11,7)
\put(-8.5,10.5){\makebox(0,0){\mbox {\scriptsize $e$}}}
\end{picture}, e
\right) :=
\left(
\unitlength.4cm
\begin{picture}(2,1)
\put(0,0){\bezier{50}(0,0)(.5,.5)(1,1)}
\put(1,0){\bezier{50}(0,1)(.5,.5)(1,0)}
\put(1,0){\line(0,1){1}}
\end{picture}\ , 1
\right)
\mbox {\hskip 1mm and \hskip 1mm}
p\left(
\unitlength 1mm
\linethickness{0.4pt}
\begin{picture}(10,4)(7,7)
\put(7,7){\line(1,1){4}}
\put(11,11){\line(1,-1){4}}
\put(15,7){\line(0,1){0}}
\qbezier(9,9)(10,8)(11,7)
\put(8.5,10.5){\makebox(0,0)[lb]{\mbox {\scriptsize $e$}}}
\end{picture}, e
\right) := 0,
\]
where $e$ denotes a metric edge of $\normalW_3$.
Finally, for $n=4$,
\begin{eqnarray*}
\label{A0}
p\left(
\unitlength 1mm
\linethickness{0.4pt}
\begin{picture}(9.7,4.5)(2.3,6)
\put(3,6){\line(0,1){0}}
\put(3,6){\line(1,1){4}}
\put(7,10){\line(1,-1){4}}
\put(7,10){\line(-1,-2){2}}
\put(5,6){\line(0,1){0}}
\put(7,10){\line(0,1){0}}
\put(7,10){\line(1,-2){2}}
\end{picture}, 1
\right)
&:=&
\left(
\unitlength 1mm
\linethickness{0.4pt}
\begin{picture}(15,6)(2,3.5)
\put(3,3){\line(0,1){0}}
\put(3,3){\line(1,1){6}}
\put(9,9){\line(1,-1){6}}
\qbezier(5,5)(6,4)(7,3)
\qbezier(7,7)(9,5)(11,3)
\put(5,7){\makebox(0,0)[cc]{\mbox {\scriptsize $b$}}}
\put(7,9){\makebox(0,0)[cc]{\mbox {\scriptsize $a$}}}
\end{picture}, a \land b
\right),
\\
\label{A1}
p\left(
\unitlength 1.4mm
\linethickness{0.4pt}
\begin{picture}(9,4)(3,4)
\put(3,4){\line(1,1){4}}
\put(7,8){\line(1,-1){4}}
\put(9,6){\bezier{50}(0,0)(-1,-1)(-2,-2)}
\put(9,6){\line(0,-1){2}}
\put(9.1,7.1){\makebox(0,0){\mbox{\scriptsize $e$}}}
\end{picture},e
\right)
&:=&
\left(
\unitlength 1mm
\linethickness{0.4pt}
\begin{picture}(13,6)(3.5,3)
\put(3,2){\line(1,1){6}}
\put(9,8){\line(1,-1){6}}
\put(9,8){\line(0,-1){4}}
\put(9,4){\line(0,1){0}}
\qbezier(9,4)(8,3)(7,2)
\qbezier(9,4)(10,3)(11,2)
\put(10.1,5){\makebox(0,0)[cc]{\mbox{\scriptsize $a$}}}
\end{picture},a
\right)
+
\left(
\unitlength 1.4mm
\linethickness{0.4pt}
\begin{picture}(9,4)(-11,5)
\put(-3,4){\line(-1,1){4}}
\put(-7,8){\line(-1,-1){4}}
\put(-9,6){\bezier{50}(0,0)(1,-1)(2,-2)}
\put(-9,6){\line(0,-1){2}}
\put(-9.1,7){\makebox(0,0){\mbox{\scriptsize $a$}}}
\end{picture},a
\right),
\\
\label{B1}
p\left(
\unitlength 1mm
\linethickness{0.4pt}
\begin{picture}(13,6)(3.5,3)
\put(3,2){\line(1,1){6}}
\put(9,8){\line(1,-1){6}}
\put(9,8){\line(0,-1){4}}
\put(9,4){\line(0,1){0}}
\qbezier(9,4)(8,3)(7,2)
\qbezier(9,4)(10,3)(11,2)
\put(10,5){\makebox(0,0)[cc]{\mbox{\scriptsize $e$}}}
\end{picture},e
\right)
&:=&
\left(
\unitlength 1.4mm
\linethickness{0.4pt}
\begin{picture}(9,4)(-11,5)
\put(-3,4){\line(-1,1){4}}
\put(-7,8){\line(-1,-1){4}}
\put(-9,6){\bezier{50}(0,0)(1,-1)(2,-2)}
\put(-9,6){\line(0,-1){2}}
\put(-9.1,7){\makebox(0,0){\mbox{\scriptsize $a$}}}
\end{picture},a
\right),
\\
\label{C1}
p\left(
\unitlength .9mm
\linethickness{0.4pt}
\begin{picture}(15,6)(2.5,3)
\put(3,2){\line(1,1){6}}
\put(9,8){\line(1,-1){6}}
\put(15,2){\line(0,1){0}}
\put(9,8){\line(0,-1){6}}
\put(9,2){\line(0,1){0}}
\qbezier(13,4)(12,3)(11,2)
\put(12,7){\makebox(0,0)[cc]{\makebox{\scriptsize $e$}}}
\end{picture}, e
\right)
&:=&
\left(
\unitlength .9mm
\linethickness{0.4pt}
\begin{picture}(15,6)(-16,3)
\put(-3,2){\line(-1,1){6}}
\put(-9,8){\line(-1,-1){6}}
\put(-9,8){\line(0,-1){6}}
\qbezier(-13,4)(-12,3)(-11,2)
\put(-12,7){\makebox(0,0)[cc]{\makebox{\scriptsize $a$}}}
\end{picture}, - a
\right)
\mbox {\hskip 2mm and}
\\
\label{A2}
p\left(
\unitlength 1mm
\linethickness{0.4pt}
\begin{picture}(15,6)(-15.5,3.5)
\put(-3,3){\line(-1,1){6}}
\put(-9,9){\line(-1,-1){6}}
\qbezier(-5,5)(-6,4)(-7,3)
\qbezier(-7,7)(-9,5)(-11,3)
\put(-5,7){\makebox(0,0){\mbox {\hskip .5mm \scriptsize $f$}}}
\put(-7,9){\makebox(0,0){\mbox {\scriptsize $e$}}}
\end{picture}, e \land f
\right)
&:=&
\left(
\unitlength 1mm
\linethickness{0.4pt}
\begin{picture}(9.7,4.5)(2.3,6)
\put(3,6){\line(0,1){0}}
\put(3,6){\line(1,1){4}}
\put(7,10){\line(1,-1){4}}
\put(7,10){\line(-1,-2){2}}
\put(5,6){\line(0,1){0}}
\put(7,10){\line(0,1){0}}
\put(7,10){\line(1,-2){2}}
\end{picture}, 1
\right),
\end{eqnarray*}
where $e$ and $f$ are metric edges of $\normalW_4$. The above
equations can be written in a more condensed form as
\begin{eqnarray*}
&
p(
\hskip .7mm
\unitlength1pt\begin{picture}(6,12)(-1.5,1)
\put(0,0){\line(0,1){10}}
\end{picture}) = \unitlength1pt\begin{picture}(6,12)(-1.5,1)
\put(0,0){\line(0,1){10}}
\end{picture},\
p(\gen 12) = \gen 12,\
p(\gen 13) = \ZbbZb,\
p(\bZbbZ) = \gen 13,\
p(\ZbbZb) = 0,&
\\
&
p(\gen 14) = \ZbbZbZb,\
p(\bZbbbZ) = \bZbbZb + \ZbbbZb,\
p(\bZbbZb) = \ZbbbZb,\
p(\bbZbbZ) = - \ZbbZbb \mbox { and }
p(\bZbZbbZ) = \gen 14,
&
\end{eqnarray*}
with the convention that binary trees are endowed with their canonical
orientations, corollas are oriented with the $+$ sign and trees $T$
with one binary and one ternary vertex are oriented as $(T,e)$, where
$e$ denotes the unique interior edge of $T$.  }
\end{example}

Let us close this section by the following proposition whose proof is
postponed to Section~\ref{postponed}.

\begin{proposition}
\label{pisu_v_susarne}
Let $(T,\omega_T)$ be an oriented fully metric tree, then
\begin{equation}\label{P}
p(\pa_{\normalW} (T,\omega_T)=\pa_{\normalK}p(T,\omega_T).
\end{equation}
Since $p$ is a operad morphism, this implies that $p$ commutes
with the differential on $\cc_*(\normalW_*)$ and therefore is a morphism
of operads in the category of chain complexes.
\end{proposition}


\section{The Saneblidze-Umble diagonal}
\label{ElisAbeth}

In this section we define the SU-diagonal~\cite{umble-saneblidze}.  Let us
start with a definition of the cubical diagonal $\Delta_\normalW$
adapted from~\cite[Section~2]{serre:AM51}:
\begin{eqnarray}
\label{pristi_utery_prohlidka}
\lefteqn{
\Delta_\normalW(T,e_1 \land \squeezedcdots \land e_k) := }
\\
\nonumber
&&:= \sum_{L,R} (-1)^{\rho_{L,R}} (T/e_L,e_1 \land
\squeezedcdots \hat{e}_{i_1} \squeezedcdots
\hat{e}_{i_l} \squeezedcdots \land e_k) \ot (T_R, e_1 \land
\squeezedcdots \hat{e}_{j_1} \squeezedcdots
\hat{e}_{j_r} \squeezedcdots \land e_k),
\end{eqnarray}
where the summation runs over all disjoint decompositions
$L \sqcup R = \{\Rada i1l\} \sqcup \{\Rada j1r\}$ of $\{\rada 1k\}$ into
ordered subsets, $T/e_L$ is the tree obtained from $T$ by contracting
edges $\{e_i;\ i \in L\}$, $T_R$ is the tree obtained by
changing the metric edges $\{e_j;\ j \in R\}$ to non-metric ones, and
$\rho_{L,R}$ is the number of couples $i \in L$, $j \in R$ such that $i < j$.
We leave as an exercise to prove:

\begin{proposition}
  The diagonal~(\ref{pristi_utery_prohlidka}) is co-associative and
  commutes with the $\circ_i$-operations introduced in~(\ref{circicw}),
  therefore
  the $W$-construction $(\ucw,\Delta_\normalW)$ is a Hopf non-$\Sigma$
  operad.
\end{proposition}

The SU-diagonal is then defined by formula~(\ref{zitra_prohlidka}),
that is
\begin{equation}
\label{AA}
\Deltasu :=(p\ot p)\circ\Delta_{\normalW}\circ q.
\end{equation}

\begin{exercise}
{\rm
Derive from definition that, in the shorthand of
Example~\ref{minuly_tyden_jsem_byl_s_Jitkou},
\begin{eqnarray}
\Deltasu(\gen 12) &=& \gen 12 \ot \gen 12, \nonumber
\\
\label{ElizabetH}
\Deltasu(\gen 13) &=& \ZbbZb \ot \gen 13 + \gen 13 \ot \bZbbZ,
\\
\nonumber
\Deltasu(\gen 14) &=& \ZbbZbZb \ot \gen 14 + \ZbbbZb \ot \bZbbZb
+ \ZbbbZb \ot \bZbbbZ + \bZbbZb \ot \bZbbbZ - \ZbbZbb \ot \bbZbbZ
+ \gen 14 \ot \bZbZbbZ, \mbox { etc.}
\end{eqnarray}
Prove also that $\Deltasu(c(n),1)$ always contains the terms
\[
(\min(n),\xi_{\min(n)})\otimes (c(n),1)\ \mbox { and }\ (c(n),1)\otimes
(\max(n),\xi_{\max(n)}).
\]
}
\end{exercise}

Let us analyze formula~(\ref{AA}) applied to $(c(n),1)$. The map
$q_n$ applied to the oriented corolla $(c(n),1) \in
\cc_{n-2}(\normalK_n)$ is, by definition, the sum of all fully
metric binary trees with standard orientations. The diagonal
$\Delta_{\normalW_n}$ acts on such a tree $(S,\omega_S)$ as
follows. Divide interior edges of $S$ into two disjoint groups,
$\{\Rada f1s\}$, $\{\Rada e1t \}$, $t+s = n-2$, and let
\[
\omega_S = \eta \cdot e_1 \land \cdots \land e_t \land f_1 \land
\cdots \land f_s,
\]
with some $\eta \in \{-1,1\}$, be the standard orientation.

Then $\Delta_{\normalW_n}(S)$ contains the term
$(S_L,\landRada e1t) \ot (S_R,\landRada f1s)$, where $S_L = S/\{\Rada f1s\}$
and $S_R$ is obtained by replacing edges $\{\Rada e1t \}$ of $S$ by
non-metric ones. We must then evaluate
\begin{equation}
  \label{eq:1}
  p_n(S_L,\landRada e1t) \ot p_n(S_R,\landRada f1s).
\end{equation}
One can also describe the pair $S_L,S_R$ as follows:  $S_R$ is the
same binary tree as  $S$, but with  only a subset of the edges
retaining the  metric label, $S_L$ is the fully-metric tree formed
from $S$ by collapsing the same subset of edges.

Let us pause a little and observe that the expression in~(\ref{eq:1})
is nonzero only for trees $S$ of a very special form.  Since the value
$p(U,\omega)$ is, for a binary fully metric tree $U$, nonzero only
when $U$ is maximal, $p_n(S_R,\landRada f1s)$ is nontrivial only when
$S_R$ is build from maximal binary fully metric trees, using the
$\circ$-operation $t$-times.  Similarly, as we saw in
Section~\ref{maps}, $p_n(S_L,\landRada e1t)$ is nonzero if and only if
$\minn{(S_L)}$ has exactly $t$ left leaning edges.

A moment's reflection convinces us that the above two conditions are
satisfied if and only if $S_R$ is build up from $t+1$ fully metric
maximal binary trees, using $t$ times $\circ_i$-operations with $i
\geq 2$ (that is, $\circ_1$ is forbidden).  Clearly $p_n(S_R,\landRada
f1s)$ is then an $n$-tree created from $t+1$ corollas using $\circ_i$
with $i \geq 2$.  Let $M^t_n$ denote the set of such $n$-trees and
$M_n := M_n^0 \sqcup \cdots \sqcup M_n^{n-2}$. A more formal
definition is that $T \in M_n$ if and only if $T = c(k) \circ_i S$ for
some $S \in M_l$, where $k+l = n-1$ and $i \geq 2$.  For example,
\[
M_1 = \{\hskip .7mm
\unitlength1pt\begin{picture}(6,12)(-1.5,1)
\put(0,0){\line(0,1){10}}
\end{picture} \},\
M_2 = \{\gen 12\},\
M_3 = \{\bZbbZ, \gen 13 \},\
M_4 = \{ \bZbZbbZ, \bZbbZb, \bZbbbZ, \bbZbbZ, \gen 14\},\ \mbox { etc.}
\]
We recommend to prove as an exercise that $M_n$ is the set of all $n$
trees $T$ whose number of interior edges is the same as the number of
left leaning edges of $\minn{T}$.

Let us reverse the process and start with an oriented $n$-tree
$(T,\xi) \in \cc_s(K_n)$ such that $T \in M^t_n$ and $\xi =
\landRada e1t$. Let $\widetilde T$ be the tree obtained from $T$
by filling all non-binary vertices by left-leaning metric edges.
Let us denote these newly created metric edges $\Rada f1s$.
Observe that
\[
p_n(\widetilde T_R, \landRada f1s) = \epsilon \cdot (T,\xi)
\]
for some $\epsilon \in \{-1,+1\}$. Define $\eta_T \in \{-1,+1\}$ by
demanding $ \eta_T \cdot \epsilon \cdot \landRada e1t \land \landRada
f1s $ to be the standard orientation of $\widetilde T$. It is not hard
to prove that $\eta_T$ indeed depends only on $T$ and not on the
choices of the labels $\Rada e1t,\Rada f1s$, as suggested by the
notation. For example, $\eta_T = 1$ for all trees from $T \in M_n$ with $n
\leq 1$ except $T=\bZbZbbZ = b(4)$ for which $\eta_T = -1$. More
generally, $\eta_{\max(n)} = (-1)^{(n-2)(n-3)/2}$.

Observe finally that $\widetilde T_L = T$. Equation~(\ref{AA}) can
then be rewritten as
\begin{equation}
\label{uz_je_tady_podzim}
\Deltasu(c(n),1)
= \sum_\doubless{T \in M^t_n}{0 \leq t \leq n-1}
  \eta_T \cdot p_n(T,\landRada e1t) \ot (T,\landRada e1t).
\end{equation}
Let us notice that the above display contains the symbol $(T,\landRada
e1t)$ twice. The first occurrence of this symbol denotes a cell of
$\cc_t(W_n)$, the second occurrence a cell of $\cc_s(K_n)$. The sign
$\eta_T$ then accounts for the difference between these two
interpretations of the same symbol.

We already observed in Example~\ref{osculetur_me} that, modulo
orientations, $p_n(T,\landRada e1t)$ in~(\ref{uz_je_tady_podzim}) is
the sum of all $n$-trees $U$ with $s$ interior edges such that
$\maxx U \leq \minn T$. This leads to the following formula for the
SU-diagonal whose spirit is closer to~\cite{umble-saneblidze}:
\begin{equation}
\Deltasu(c(n),1) = \sum \vartheta \cdot (U,\omega_U) \ot (T,\omega_T),
\end{equation}
where, as usual, $c(n)$ is the $n$-corolla representing the top
dimensional cell of $K_n$, the summation is taken over
all $(U,\omega_U)$, $(T,\omega_T)$
with $\maxx U \leq \minn T$ and $\dim(S,\omega_S) + \dim(T,\omega_T) = n$,
and $\vartheta$ is a sign which can be picked up by comparing this
formula to~(\ref{uz_je_tady_podzim}).

\section{Non-existence of a co-associative diagonal}
\label{tyden-pred-prohlidkou}

As we already indicated, the SU-diagonal is not co-associative, that is,
\begin{equation}\nonumber
(\Deltasu \ot \id)\Deltasu \not = (\id \ot \Deltasu)\Deltasu.
\end{equation}
While still
\[
(\Deltasu \ot \id)\Deltasu(\gen12) = (\id \ot \Deltasu)\Deltasu(\gen12)
\mbox { and }
(\Deltasu \ot \id)\Deltasu(\gen13) = (\id \ot \Deltasu)\Deltasu(\gen13),
\]
the co-associativity breaks already for $\gen 14$, explicitly:
\begin{equation}
\label{elisAbeth}
(\Deltasu \ot \id)\Deltasu(\gen14) - (\id \ot
\Deltasu)\Deltasu(\gen14)
= \pa (\ZbbbZb \ot \bZbbZb \ot \bZbbbZ)
\end{equation}
The SU diagonal is also not co-commutative. This means that
\[
T(\Deltasu) \not= \Deltasu,
\]
where $T : \opinf \ot \opinf \to \opinf \ot \opinf$ is the `flip.'
More explicitly, while $T(\Deltasu)(\gen 12) =   \Deltasu(\gen 12)$,
\[
\Deltasu(\gen 13)-  T(\Deltasu)(\gen 13) = \pa (\gen 13 \ot \gen 13).
\]

In the rest of this section we show that the non-coassociativity of
$\Deltasu$ is not due to bad choices in the definition, but follows
from a deeper principle, namely:

\begin{theorem}
\label{ElisaBeth}
The operad $\opinf$ does not admit a co-associative diagonal.
Therefore the operad $\opinf$ for \ainfty-algebras is not a Hopf
operad in the sense of~\cite{getzler-jones:preprint}.
\end{theorem}

\noindent
{\bf Proof.}
The proof is boring and the reader is warmly encouraged to skip it.
The idea is to try to construct inductively a
co-associative diagonal $\Delta$ and observe that at a certain stage there is a
non-trivial co-associativity constraint. Let us start with the construction.
For $\gen 12$ we are forced to take
\[
\Delta(\gen 12) := \gen 12 \ot \gen 12.
\]
The most general form of $\Delta(\gen 13)$ is
\[
\Delta(\gen 13) = (a \ZbbZb + b \bZbbZ) \ot \gen 13 +
                  \gen 13 \ot (c \ZbbZb + d \bZbbZ),
\]
with some $a,b,c,d \in \bfk$. The compatibility with the differential
$\pa$ of $\opinf$ means that
\[
\pa \Delta(\gen 13) = (a \ZbbZb + b \bZbbZ) \ot (\bZbbZ - \ZbbZb) +
                       (\bZbbZ - \ZbbZb) \ot  (c \ZbbZb + d \bZbbZ)
\]
must be the same as
\[
\Delta(\pa \gen 13) = \Delta(\bZbbZ - \ZbbZb)
                    =  \bZbbZ \ot\bZbbZ - \ZbbZb \ot  \ZbbZb.
\]
This is clearly equivalent  to
\[
a+c =1,\ b + d =1,\  a=d \mbox { and } b =c.
\]
It can be equally easily verified that the co-associativity
\[
(\Delta \ot \id)\Delta(\gen 13) = (\id  \ot \Delta)\Delta(\gen 13)
\]
is equivalent to
\[
a = a^2,\ b = b^2,\ c = c^2,\ d = d^2,\ ab = 0 \mbox { and } cd = 0.
\]
We conclude that the only two co-associative solutions
are either $(a,b,c,d) = (1,0,0,1)$ or $(a,b,c,d) = (0,1,1,0)$, that is
either
\begin{equation}
\label{Elizabeth}
\Delta (\gen 13) = \ZbbZb \ot \gen 13 + \gen 13 \ot \bZbbZ.
\end{equation}
or
\begin{equation}
\label{ELIzabeth}
\Delta (\gen 13) = \bZbbZ \ot \gen 13 + \gen 13 \ot \ZbbZb.
\end{equation}
Let us assume solution~(\ref{Elizabeth}) which coincides with the SU-diagonal
(compare~(\ref{ElizabetH})) -- solution~(\ref{ELIzabeth}) is just the flip
$T(\Deltasu(\gen 13))$ and this case can be discussed by flipping all the
steps below.
We will be looking for $\Delta$ of the form
$\Delta = \Deltasu + \delta$ with some perturbation $\delta : \opinf \to \opinf
\ot \opinf$ satisfying, of course, $\delta(\gen 12) = \delta(\gen 13)
= 0$. Since we know that $\Deltasu$ is a chain map, $\delta$ must be
a chain map as well.

Observe that $\delta(\gen 14)$ depends on $35$ parameters. Therefore
the co-associativity of $\Delta$ and the chain condition on $\delta$
is expressed by a system of linear equations in $35$ variables! We are
going to show that this system has no solution. This might be a
formidable task, but we will simplify it by making some wise guesses.
Let us write\def\jejednadva{{J_{(1)}^i \ot J_{(2)}^i}}
\begin{equation}
\label{elisabeTH}
\delta(\gen 14) = A \ot \gen 14 + \sum_i \jejednadva + \gen 14 \ot B,
\end{equation}
where $A,B \in \opinf_0(4)$ and $\jejednadva \in \opinf_1(4) \ot \opinf_1(4)$.
Let us also denote%
\def\RHS{{\it RHS\/}}\def\LHS{{\it LHS\/}}
\[
\LHS{} : = [(\delta \ot \id)\delta + (\Deltasu \ot \id)\delta +
   (\delta \ot \id)\Deltasu + (\Deltasu \ot \id)\Deltasu \hskip .5mm](\gen 14)
\]
and
\[
\RHS{} : = [(\id \ot\delta)\delta + (\id \ot\Deltasu)\delta +
   (\id \ot \delta)\Deltasu + (\id \ot\Deltasu)\Deltasu \hskip .5mm](\gen 14).
\]
The co-associativity of $\Delta$ at $\gen 14$ of course means that $\LHS{} =
\RHS$. An easy calculation shows that the
only term of \LHS\ of the form $\gen 14 \ot \mbox
{\it something\/}$ is
\[
\gen 14 \ot (B \ot B + \bZbZbbZ \ot B + B \ot \bZbZbbZ + \bZbZbbZ \ot \bZbZbbZ)
\]
while the only term of \RHS\ of the same form is
\[
\gen 14 \ot (B \ot B + \bZbZbbZ \ot \bZbZbbZ).
\]
Associativity $\RHS=\LHS$ then evidently means
\[
\bZbZbbZ \ot B + B \ot \bZbZbbZ =0,
\]
which, since ${\it char}(\bfk) \not= 2$, clearly implies
$B=0$. Using the same trick we see also that $A=0$,
therefore $\delta(\gen 14)$ must be of the
form
\[
\delta(\gen 14) = \sum_i \jejednadva.
\]
Since $\delta$ is a chain map, trivial on $\gen 12$ and $\gen 13$,
$\pa \delta(\gen 14) = 0$, which means that
\[
0 = \pa \delta(\gen 14)= \sum_i \pa J^i_{(1)} \ot J^i_{(2)} -
\sum_i J^i_{(1)} \ot \pa J^i_{(2)}.
\]
Looking separately at the components of bidegrees $(1,0)$ and
$(0,1)$, respectively, and assuming, without loss of generality, that
the elements $J^i_{(1)}$ (resp.~$J^i_{(2)}$) are linearly independent,
we conclude that $\pa J^i_{(1)} = \pa J^i_{(2)} = 0$. Because
each cycle in $\opinf_1(4)$ is a scalar multiple
of\def\pactyri{\pa(\gen 14)} $\pactyri$, we see that
\[
\delta(\gen 14) = \alpha (\pactyri \ot \pactyri)
\]
for some scalar $\alpha \in \bfk$. So we managed to cut $35$
parameters in~(\ref{elisabeTH}) to one!
Now
\[
\LHS{} = \alpha\{
\Deltasu(\pactyri) \ot \pactyri + \pactyri \ot \pactyri \ot \bZbZbbZ\} +
(\Deltasu \ot \id)\Deltasu (\gen 14)
\]
and
\[
\RHS{} = \alpha\{
\pactyri \ot \Deltasu(\pactyri) + \ZbbZbZb \ot \pactyri \ot \pactyri\}
+(\id \ot \Deltasu)\Deltasu(\gen 14).
\]
The only terms of the \LHS\ of the form $\bZbbbZ \ot \mbox {\it
something\/}$ are
\[
\alpha\{
\bZbbbZ \ot \bZbZbbZ \ot \pactyri + \bZbbbZ \ot \pactyri \ot \bZbZbbZ
\}
\]
while in the \RHS, there is only one term of this form, namely
\[
\alpha (\bZbbbZ \ot \Deltasu(\pactyri)).
\]
The only term of the form $\bZbbbZ \ot \bZbZbbZ  \ot \mbox {\it
something\/}$ in the above two displays is
\[
\alpha(\bZbbbZ \ot \bZbZbbZ \ot \pactyri)
\]
coming from the first term of the first display.
This implies that $\alpha = 0$, therefore $\delta =
0$ and $\Delta = \Deltasu$. But this is not possible, because the
co-associativity of $\Deltasu$ is violated already on $\gen 14$,
as we saw in~(\ref{elisAbeth}).%
\qed

\section{Remaining proofs}
\label{postponed}

In this section we prove Propositions~\ref{ale_porad_na_ni_myslim}
and~\ref{pisu_v_susarne}. Let us start with a

\noindent
{\bf Proof of Proposition~\ref{ale_porad_na_ni_myslim}.}
By definition of $\pa_\normalK$,
\begin{eqnarray}
\pa_\normalK(c(n),1)&=&\sum_{\doubless{r+s=n+1}{1\leq i\leq r}}
(c(r)\circ_i c(s), e)
=\sum_{\doubless{r+s=n+1}{1\leq i\leq r}}
(-1)^{r(s-2)+i(s+1)}(c(r),1)\circ_i (c(s),1)\nonumber
\\
&=&\sum_{\doubless{r+s=n+1}{1\leq i\leq r}}
(-1)^{(r+i)s +i}(c(r),1)\circ_i (c(s),1).
\nonumber
\end{eqnarray}
The sign comes from formula~(\ref{circick}) setting $l=s-2$, since
$c(s)\in \cc_{s-2}(\normalK_s)$. Applying $q$ gives
\begin{equation}
\label{A}
q(\pa_\normalK(c(n),1))=\sum_{\doubless{r+s=n+1}{1\leq i\leq r}}
(-1)^{(r+i)s +i}q(c(r),1)\circ_i q(c(s),1).
\end{equation}
The expression on the right is a sum over all binary rooted planar metric trees
 with $n$ leaves and  one non-metric edge. On the other hand,
$$\pa_\normalW(q(c(n),1))=\sum_{T\in\mbn} \pa_\normalW(T,\omega_T).$$
According to (\ref{pacw}), the terms  in  $\pa_\normalW(T,\omega_T)$ are of
two types:

\vskip 1mm
\noindent
\halign{
\vtop{\parindent=0pt\hsize=1.8cm\hfill\strut#\strut}&%
\hskip 2mm\vtop{\parindent=0pt\hsize=14cm\strut#\strut}\cr
Type A,& in which a metric edge has been changed
to a non-metric edge and%
\cr
Type B,& in which a metric edge has been collapsed,
creating a fully metric tree which is binary except for one tertiary
vertex.%
\cr
}

In the sum of type B terms the same cell appears twice with
opposite signs, since there are exactly two binary trees which
give rise to the same  tree with a unique tertiary vertex.
The terms of type A, with one non-metric edge, run over the set of all
binary rooted planar metric trees with one non-metric edge, which is
the same set as that appearing in the sum on the right of equation~(\ref{A}).
It only remains to compare the orientations of the corresponding
terms on the two sides of (\ref{qpa}).
According to Definition~\ref{defq},
\begin{eqnarray}
\label{brcircbs1}
\lefteqn{
(-1)^{(r+i)s +i}q(c(r),1)\circ_i q(c(s),1)=\hskip .5cm}
\\
&& \hskip .5cm
(-1)^{(r+i)s +i}
(\max(r)\circ_i \max(s),e_1\wedge\cdots_\wedge e_{r-2}\wedge
f_1\wedge\cdots\wedge f_{s-2})+\cdots,
\nonumber
\end{eqnarray}
where the term shown explicitly on the right is the leading order term
relative to the order relation on binary trees, $e_1,\ldots,e_{r-2}$
label the edges in $\max(r)$ and $f_1,\ldots,f_{s-2}$ label the edges
in $\max(s)$.  Since the definition of the standard orientation on an
arbitrary fully metric binary tree involves the same associativities
independent of the size of the tree, it is sufficient to compare the
orientation of the leading order term in (\ref{brcircbs1}) with the
orientation of the corresponding term in $\pa_\normalW(q(c(n),1))$.
If these orientations agree, so will the orientations of all the other
terms.

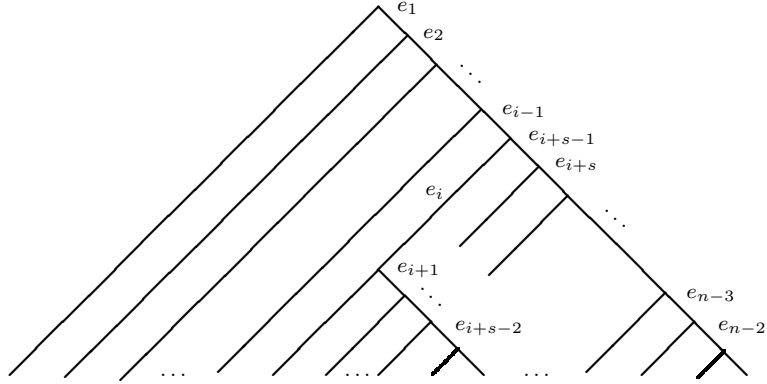
\begin{figure}
\begin{center}
\unitlength 3.5mm
\thicklines
\begin{picture}(29,15)(.5,0)
\put(1,1){\line(1,1){14}}
\put(15,15){\line(1,-1){14}}
\put(16.1,13.9){\line(-1,-1){13}}
\put(17.2,12.8){\line(-1,-1){12}}
\put(18.9,11.1){\line(-1,-1){10}}
\put(20,10){\line(-1,-1){9}}
\put(15,5){\line(1,-1){4}}
\put(27,3){\line(-1,-1){2}}
\put(25.9,4.1){\line(-1,-1){3}}
\put(16,4){\line(-1,-1){3}}
\put(18,2){\bezier{50}(-1,-1)(-.5,-.5)(0,0)}
\put(17,3){\line(-1,-1){2}}
\put(21.1,8.9){\line(-1,-1){3}}
\put(22.2,7.8){\line(-1,-1){3}}
\put(28.1,1.9){\bezier{50}(-1,-1)(-.5,-.5)(0,0)}
\put(15.7,14.7){\makebox(0,0)[lb]{\scriptsize $e_1$}}
\put(16.7,13.7){\makebox(0,0)[lb]{\scriptsize $e_2$}}
\put(18,12){\makebox(0,0)[lb]{\scriptsize $\ddots$}}
\put(19.7,10.7){\makebox(0,0)[lb]{\scriptsize $e_{i-1}$}}
\put(17.5,7.7){\makebox(0,0)[rb]{\scriptsize $e_i$}}
\put(15.7,4.7){\makebox(0,0)[lb]{\scriptsize $e_{i+1}$}}
\put(16.5,3.5){\makebox(0,0)[lb]{\scriptsize $\ddots$}}
\put(17.9,2.5){\makebox(0,0)[lb]{\scriptsize $e_{i+s-2}$}}
\put(20.7,9.7){\makebox(0,0)[lb]{\scriptsize $e_{i+s-1}$}}
\put(21.7,8.7){\makebox(0,0)[lb]{\scriptsize $e_{i+s}$}}
\put(23.5,6.5){\makebox(0,0)[lb]{\scriptsize $\ddots$}}
\put(26.7,3.7){\makebox(0,0)[lb]{\scriptsize $e_{n-3}$}}
\put(27.9,2.5){\makebox(0,0)[lb]{\scriptsize $e_{n-2}$}}
\put(7.2,1){\makebox(0,0){\scriptsize $\ldots$}}
\put(14.2,1){\makebox(0,0){\scriptsize $\ldots$}}
\put(21,1){\makebox(0,0){\scriptsize $\ldots$}}
\end{picture}
\caption{%
The binary tree in the figure is derived from the maximal binary tree
$\max(n)$ by moving the $s-1$  adjacent vertices between edges $e_i$ and
$e_{i+s-1}$; therefore, its standard orientation is
$(-1)^{s-1} e_1\wedge\cdots\wedge e_{n-2}$.}
\label{brcircbs}
\end{center}
\end{figure}

Assume $i < r$. Applying $\pa_\normalW$ to the fully metric binary
tree with standard orientation appearing in Figure~\ref{brcircbs}, we
get (among others) the term $(-1)^{(i+s)}(\max(r)\circ_i
\max(s),e_1\wedge\cdots\hat{e}_i \ldots\wedge e_{n-2})$ with $e_i$
changed to a non-metric edge, and the edges labeled $e_{i+1},\ldots,
e_{i+s-2}$ corresponding to the edges in $\max(s)$.  Reordering the terms
in the orientation element appearing in (\ref{brcircbs1}) so that
$f_1,\ldots, f_{s-2}$ appear in sequence between $e_{i-1}$ and $e_i$
introduces a sign factor $(-1)^{(s-2)(r-i-1)}$. But
$(-1)^{(s-2)(r-i-1)}(-1)^{(r+i)s +i} =(-1)^{i+s}$, since so the signs
agree.
For $i=r$, when $\max(r)\circ_r \max(s) = \max(r+s-1)$, the analysis is much
easier and we leave it to the reader.%
\qed

\noindent
{\bf Proof of Proposition~\ref{pisu_v_susarne}.}  The case $n=2$ is
trivial. Assuming the proposition is true for fully metric trees
$(T,\omega_T)\in \cc_*(\normalW_m)$ for $m<n$, we will prove it for
$\cc_k(\normalW_n)$, starting with $k=n-2$ and descending.  In the
case $\cc_{n-2}(\normalW_n)$, which involves binary fully metric
trees, we begin with the maximal binary metric tree.  We need to prove
the commutativity of Figure~\ref{pcommute}, which follows from the
equations in Figure~\ref{pcommute2} once we check the signs.
\begin{figure}
\unitlength 5mm
\thicklines
\begin{picture}(25,12)(0,0)
\put(0,8.9){
\unitlength 1.8mm
\put(0,1){
\put(3,2){\line(1,1){6}}
\put(9,8){\line(1,-1){6}}
\put(10.5,6.5){\line(-1,-1){4.5}}
\put(13,4){\line(-1,-1){2}}
\put(11.5,6){\makebox(0,0)[rt]{\scriptsize$\ddots$}}
}
\put(2,6){\makebox(0,0)[cr]{$\leftrubber8$}}
\put(16,6){\makebox(0,0)[lc]{$,e_1 \land \cdots \land e_{n-2}\rightrubber8$}}
}
%
\put(23,7.9){
\unitlength 2.4mm 
\linethickness{0.4pt}
\begin{picture}(14,9)(0,0)
\put(6,5){\line(1,1){4}}
\put(10,9){\line(0,1){0}}
\put(10,9){\line(1,-1){4}}
\put(14,5){\line(0,1){0}}
\put(10,9){\line(-1,-2){2}}
\put(10,9){\line(1,-2){2}}
\put(10,5){\makebox(0,0)[cc]{\scriptsize $\cdots$}}
\end{picture}
\put(-8,7){\makebox(0,0)[cr]{$\leftrubber8$}}
\put(1,7){\makebox(0,0)[lc]{$,1\rightrubber8$}}
}
\put(2.6,-1.5){
\unitlength 1.8mm
\thicklines
\begin{picture}(31,11)(0,0)
\put(13,0){\levystromecek}
\put(-15,0){\pravystromecek}
\put(18.5,8){\makebox(0,0)[cc]{$-$}}
\put(-8.5,8){\makebox(0,0)[lc]{$\leftrubber{10}\sum (-1)^i \leftrubber6$}}
\put(31.5,8){\makebox(0,0)[lc]
{$\rightrubber6,e_1\land\squeezedcdots \hskip .87mm {\hat e}_i
\squeezedcdots\land e_{n-2}\rightrubber{10}$}}
\end{picture}
}
\put(28.5,.4){\makebox(0,0)[cc]{
\unitlength 1.8mm
\thicklines
\begin{picture}(18,11)(0,0)
\put(4,4){\line(1,1){7}}
\put(11,11){\line(1,-1){7}}
\put(18,4){\line(0,1){0}}
\put(11,11){\line(-3,-4){5.25}}
\put(11,11){\line(3,-4){5.25}}
\put(11,11){\line(0,-1){4}}
\put(11,7){\line(-1,-1){3}}
\put(11,7){\line(1,-1){3}}
\put(11,7){\line(-1,-3){1}}
\put(6.3,4){\multiput(0,0)(.4,0){3}{\scriptsize $.$}}
\put(12,4.1){\makebox(0,0)[cc]{\scriptsize $\cdots$}}
\put(14.4,4){\multiput(0,0)(.4,0){3}{\scriptsize $.$}}
\put(10.6,8.5){\makebox(0,0)[cr]{\scriptsize $e$}}
\put(17,8){\makebox(0,0)[lc]{$,e \rightrubber{10}$}}
\put(4.5,8){\makebox(0,0)[rc]{$\leftrubber{10} \sum$}}
\end{picture}
}}
\put(19.5,12){\makebox(0,0)[cc]{$p$}}
\put(3,6){\makebox(0,0)[cc]{$\partial$}}
\put(24,1.9){\makebox(0,0)[cb]{$p$}}
\put(30,6){\makebox(0,0)[lc]{$\partial$}}
\put(4,8){\vector(0,-1){4}}
\put(29,8){\vector(0,-1){4}}
\put(14,11){\vector(1,0){10}}
\put(22.8,1.2){\vector(1,0){2}}
\end{picture}
\caption{\label{pcommute}%
The commutative diagrams for
the chain map $p$ in degrees $n-2$ acting on the
maximal binary metric tree.}
\end{figure}
\begin{figure}
\begin{center}
\thicklines
{
\unitlength=1.000000pt
\begin{picture}(210.00,130.00)(0.00,10.00)
\put(0,60){
\put(70,15){\unitlength=.8pt
\put(158.00,-19.00){\makebox(0.00,0.00)%
    {$\underbrace{\rule{17mm}{0mm}}_{(i+1) \mbox {\scriptsize{} leaves}}$}}
\put(172.00,5.00){\makebox(0.00,0.00)%
    {$\underbrace{\rule{8.5mm}{0mm}}_{s \mbox {\scriptsize{} leaves}}$}}
\put(6.6666,0){
\put(170.00,40.00){\line(1,-2){10.00}}
\put(170.00,40.00){\line(-1,-2){10.00}}
\put(170.00,40.00){\line(-1,-1){20.00}}}
\put(170.00,60.00){\line(1,-3){6.6666}}
\put(170.00,60.00){\line(-3,-4){30.00}}
\put(198,20.00){\makebox(0.00,0.00){\scriptsize $\cdots$}}
\put(176.50,20.00){\makebox(0.00,0.00){\scriptsize $\cdots$}}
\put(148.00,20.00){\makebox(0.00,0.00){\scriptsize $\cdots$}}
\put(171.00,45.00){\makebox(0.00,0.00)[r]{\scriptsize $e$}}
\put(170.00,60.00){\line(1,-1){40.00}}
\put(130.00,20.00){\line(1,1){40.00}}
}
\put(145.00,37.00){\makebox(0.00,0.00){$\displaystyle\sum_{1\leq s \leq i+1}$}}
\put(10.00,40.00){\makebox(0.00,0.00){$p \leftrubber{10}$}}
{\unitlength 1.8mm
\put(-15,5){\makebox(0.00,0.00){$\pravystromecek$}}
}
\put(90.00,40.00){\makebox(0.00,0.00){$\rightrubber{10}$}}
\put(120.00,40.00){\makebox(0.00,0.00){$=$}}
}
\put(0,-85){
\put(40,23){\unitlength=.8pt
\put(198,90.00){\makebox(0.00,0.00)[b]{\scriptsize $\cdots$}}
\put(170.00,90.00){\makebox(0.00,0.00)[b]{\scriptsize $\cdots$}}
\put(145.00,90.00){\makebox(0.00,0.00)[b]{\scriptsize $\cdots$}}
\put(190.00,120.00){\makebox(0.00,0.00){\scriptsize $e$}}
\put(190.00,110.00){\line(0,-1){20.00}}
\put(190.00,110.00){\line(-1,-2){10.00}}
\put(170.00,130.00){\line(-1,-4){10.00}}
\put(170.00,130.00){\line(1,-1){40.00}}
\put(130.00,90.00){\line(1,1){40.00}}
}
\put(10.00,110.00){\makebox(0.00,0.00){$p \leftrubber{10}$}}
{\unitlength 1.8mm
\put(-2,18){\makebox(0.00,0.00){$\levystromecek$}}
}
\put(90,110.00){\makebox(0.00,0.00){$\rightrubber{10}$}}
\put(120.00,110.00){\makebox(0.00,0.00){$=$}}
}
\end{picture}}
\end{center}
\caption{\label{pcommute2}%
  Each of the trees labeling the faces of $\normalK_n$ appears
  precisely once as a term in $\pa \max(n)$. The orientation elements
  (not shown in the figure) are $(-1)^{i-1}e_1\wedge
  \cdots\hat{e_i}\cdots\wedge e_{n-2}$ in the upper left, $(-1)^i
  e_1\wedge\cdots\hat{e_i}\cdots\wedge e_{n-2}$ in the lower left, and
  $e$ for both trees on the right.}
\end{figure}

Let us start with the second equation in Figure~\ref{pcommute2}.  The
tree in parentheses on the left with orientation element
$-(-1)^{i}e_1\wedge\cdots\hat{e}_i\cdots\wedge e_{n-2}$, which
corresponds to one of the terms appearing in
$\pa_\normalW(\max(n),\omega_{\max(n)})$, is equal to
$-(-1)^{i}(\max(i+1),\omega_{\max(i+1)})\circ_{i+1} (\max(n-i),
\omega_{\max(n-i)})$ and therefore its image under $p$ is
\begin{eqnarray*}
\lefteqn{
-(-1)^{i}p(\max(i+1),\omega_{\max(i+1)})
\circ_{i+1}p(\max(n-i),\omega_{\max(n-i)}) = \hskip 2cm }
\\
&&\hskip 2cm  =-(-1)^{i}(c(i+1),1)\circ_{i+1} (c(n-i),1)
\\
&&\hskip 2cm=-(-1)^{i+(i+1)(n-i-2)+ (i+1)(n-i+1)}
(c(i+1)\circ_{i+1}c(n-i),e)
\\
&&\hskip 2cm=(c(i+1)\circ_{i+1}c(n-i),e),
\end{eqnarray*}
as required.
The orientation element for  the trees  on the right side of the
first equation in Figure~\ref{pcommute2} with $s$ leaves is
\[
(-1)^{i+(n-3)(n-4)/2}e_1\wedge\cdots\hat{e}_i\cdots\wedge e_{n-2}\into
(-1)^{s-1} e_1\wedge\cdots\wedge e_{i-s+1}\wedge e
\wedge e_{i-s+2}\wedge\cdots \wedge  e_{n-2}=e.
\]
Thus $p$ commutes with $\pa$ on the maximal
binary fully metric tree.

Next we will show that $p$ commutes with $\pa$ for all binary fully
metric trees. To simplify notation, we will not indicate the
orientation element.  For a non-maximal fully metric binary tree $T$,
$p(T)=0$, because $\minn{T} = T$ has less than $n-2$ left leaning
edges.  The only binary fully metric trees for which $p(\pa T)\neq 0$
are trees of the type appearing in Figure~\ref{brcircbs} with only one
right-leaning internal edge.

Let $T^{i,s}$ be the tree in Figure~\ref{brcircbs}, and $T^{i,s}_j$
the term in $\pa T^{i,s}$ with edge $e_j$ non-metric. Then, for $i\neq
j$, $T^{i,s}_j$ is a $\circ$-composition of two fully metric binary
trees, one of which is not maximal. Since $p$ is a operad map, the
image $p(T^{i,s}_j)$ is also a $\circ$-composition, but one of the two
components is zero, since $p(T)=0$ when $T$ is fully metric binary but
not maximal.

For $j\neq i,i-1$ we also have $p(T^{i,s}/e_j)=0$, because the binary
tree $\minn{(T^{i,s})}$ has two right leaning edges. Thus the only
terms in $\pa T^{i,s}$ whose image under $p$ is not zero are
$T^{i,s}/e_{i-1},T^{i,s}/e_{i},$ and $ T^{i,s}_i$. It follows
immediately from the definition of $p$ that
\begin{equation}
\label{ei1}
p(T^{i,s}/e_i)= p(T^{i,s}_i)+p(T^{i,s}/e_{i-1}).
\end{equation}
In fact, the one term appearing in $p(T^{i,s}/e_i)$ and not appearing
in $p(T^{i,s}/e_{i-1})$ is $p(T^{i,s}_i)$.  Therefore,
\[
p(\pa T^{i,s})=p((-1)^{(i-1)} T^{i,s}/e_{i-1} +(-1)^{i-1}T^{i,s}_i
+(-1)^i T^{i,s}/e_i)=0=\pa p(T^{i,s}).
\]
This completes the proof of (\ref{P}) for $T\in\cc_{n-2}(\normalW_n)$.

\begin{figure}
\begin{center}
\thicklines
\unitlength=.9pt
\begin{picture}(410.00,120.00)(0.00,0.00)
%
\put(320.00,120.00){\makebox(0.00,0.00){%
    The subtree of both $\minn {T}$ and $\minn{(T/e_i)}$}}
\put(287.00,33.00){\makebox(0.00,0.00)[b]{\scriptsize $e_i$}}
\put(294.00,56.00){\makebox(0.00,0.00)[t]{\scriptsize $e_{i-1}$}}
\put(285.00,25.00){\line(1,-1){15.00}}
\put(295.00,35.00){\line(1,-1){25.00}}
\put(310.00,50.00){\line(-1,-1){40.00}}
\put(290.00,70.00){\line(1,-1){60.00}}
\put(300.00,80.00){\line(1,-1){70.00}}
\put(320.00,100.00){\line(1,-1){90.00}}
\put(230.00,10.00){\line(1,1){90.00}}
\put(250,10){\makebox(0.00,0.00){$\cdots$}}
\put(285,10){\makebox(0.00,0.00){$\cdots$}}
\put(335,10){\makebox(0.00,0.00){$\cdots$}}
\put(387,10){\makebox(0.00,0.00){$\cdots$}}
%
\put(90.00,120.00){\makebox(0.00,0.00){The relevant subtree of $T$}}
\put(88.00,40.00){\makebox(0.00,0.00)[r]{\scriptsize $e_i$}}
\put(97.00,58.00){\makebox(0.00,0.00)[r]{\scriptsize $e_{i-1}$}}
\put(83.00,27.00){\line(1,-1){18.00}}
\put(100.00,50.00){\line(1,-4){10.00}}
\put(100.00,50.00){\line(3,-4){30.00}}
\put(100.00,50.00){\line(-3,-4){30.00}}
\put(90.00,100.00){\line(-2,-5){36.00}}
\put(90.00,100.00){\line(3,-5){54.00}}
\put(90.00,100.00){\line(1,-5){10.00}}
\put(90.00,100.00){\line(1,-1){90.00}}
\put(0.00,10.00){\line(1,1){90.00}}
\put(30,10){\makebox(0.00,0.00){$\cdots$}}
\put(85,10){\makebox(0.00,0.00){$\cdots$}}
\put(120,10){\makebox(0.00,0.00){\scriptsize $\cdots$}}
\put(160,10){\makebox(0.00,0.00){$\cdots$}}
\end{picture}
\caption{\label{pcommute3}%
 The only  configuration of edges
in $T$ which is relevant in the calculation of $p(\pa T)$.
The configuration may occur at any vertex, not necessarily
at the root.}
\end{center}
\end{figure}
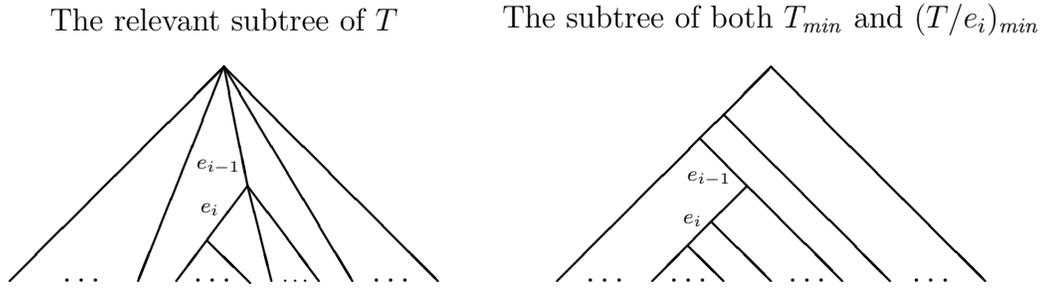

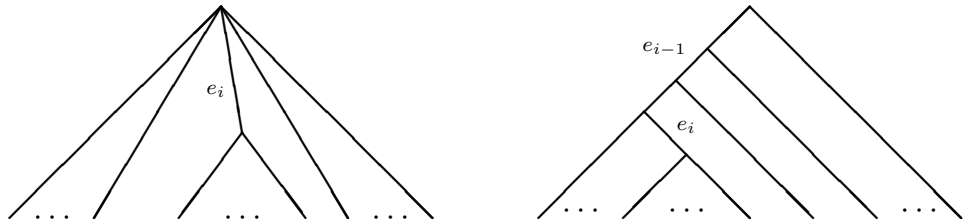
\begin{figure}
\begin{center}
{
\thicklines
\unitlength=.8pt
\begin{picture}(450.00,110.00)(0.00,0.00)
%
\put(430.00,0.00){\makebox(0.00,0.00)[b]{$\cdots$}}
\put(320.00,0.00){\makebox(0.00,0.00)[b]{$\cdots$}}
\put(270.00,0.00){\makebox(0.00,0.00)[b]{$\cdots$}}
\put(320.00,80.00){\makebox(0.00,0.00)[r]{\scriptsize $e_{i-1}$}}
\put(320.00,40.00){\makebox(0.00,0.00)[b]{\scriptsize $e_i$}}
\put(330.00,80.00){\line(1,-1){80.00}}
\put(315.00,65.00){\line(1,-1){65}}
\put(320.00,30.00){\line(-1,-1){30.00}}
\put(300.00,50.00){\line(1,-1){50.00}}
\put(350.00,100.00){\line(1,-1){100.00}}
\put(250.00,0.00){\line(1,1){100.00}}
%
\put(180.00,0.00){\makebox(0.00,0.00){$\cdots$}}
\put(110.00,0.00){\makebox(0.00,0.00){$\cdots$}}
\put(20.00,0.00){\makebox(0.00,0.00){$\cdots$}}
\put(102.00,60.00){\makebox(0.00,0.00)[r]{\scriptsize $e_i$}}
\put(100.00,100.00){\line(3,-5){60.00}}
\put(110.00,40.00){\line(3,-4){30.00}}
\put(110.00,40.00){\line(-3,-4){30.00}}
\put(100.00,100.00){\line(1,-6){10.00}}
\put(100.00,100.00){\line(-3,-5){60.00}}
\put(100.00,100.00){\line(1,-1){100.00}}
\put(0.00,0.00){\line(1,1){100.00}}
\end{picture}}
\caption{\label{pcommute4}%
  Collapsing the edge $e_{i-1}$ from the subtree on the left of
  Figure~\ref{pcommute3} creates the subtree $S$ shown here on the
  left. Filling-in to get a binary subtree $\minn{S}$ creates the
  subtree of $\minn{(T/e_{i-1})}$ shown here on the right, with
  edge $e_i$ left-leaning and $e_{i-1}$ right-leaning}
\end{center}
\end{figure}

Now, assuming that (\ref{P}) is true for all $T\in \cc_j(\normalW_n)$
for $k<j\leq n-2$ for all $T\in\cc_*(\normalW_m)$ for $m<n$, we need
to prove it for $T\in \cc_k(\normalW_n)$.  Let $T$ be a fully metric
tree with $k$ edges labeled $e_1,\ldots,e_k$ and $\minn{T}$ the binary
tree given by filling in, and label the $k$ edges in $\minn{T}$
corresponding to the original edges by the same labels.  All the other
edges of $\minn{T}$ are right leaning.  If less than $k-1$ of the
edges $e_1,\ldots, e_k$ in $\minn{T}$ are left-leaning, then
$p(T)=0=p(\pa T)$ and therefore, $\pa p(T)=p (\pa T)$.  Suppose first
that $\minn{T}$ has $k-1$ left-leaning edges, and $e_i$ is
right-leaning.  Just as for the binary metric trees, the only tree in
$\pa T$ for which the image under $p$ is non-zero are $T/e_{i-1},
T/e_i$ and $T_i$, where $e_{i-1}$ and $e_i$ are adjacent edges in $T$.
The configuration is illustrated in Figure~\ref{pcommute3}. The
subtree on the right of Figure~\ref{pcommute3} appears as a subtree in
both in $\minn{T}$ and $(T/e_i)_{\it bin}$, and the tree on the right
of Figure~\ref{pcommute4} appears as a subtree in
$\minn{(T/e_{i-1})}$. The following equation analogous to (\ref{ei1})
applies in this case
\begin{equation}
\label{ei2}
p(T/e_i)= p(T_i)+p(T/e_{i-1}).
\end{equation}
and the remainder of the proof of (\ref{P}) in this case is the same
as before.

Next we consider the case when there are exactly $k$ left-leaning
edges in $\minn{T}$. In general, for any $k$ cell $T$ such that
$\minn{T}$ has $k$ left-leaning edges we can choose a $k+1$ cell
$T^+$ such that $\pa T^+$ contains $T$ as a summand and all other
summands of the type $R_i:=T^+/e_i$ for $i=1,\ldots,k$ have the
property that $\minn{(R_i)}$ has $k-1$ left-leaning edges.  The
tree $T^+$ can be defined as follows: pick any non-binary vertex $v$
in $T$ with $r\geq 3$ incoming edges and replace the corolla with
vertex $v$ by subtree of type $c(2)\circ_1 c(r-1)$ with the new edge
labeled $e_{k+1}$.  Then $T^+/e_{k+1}=T$ and for $i=1,\ldots,k$,
$\minn{(T^+/e_i)}$ has $k-1$ left-leaning edges, since the edge
corresponding to $e_{k+1}$ in $\minn{(T^+/e_i)}$ is right-leaning.
Denote the faces of type $1$ $T^+_j$ in which a metric edge is changed
to a non-metric edge by $S_j$, $j=1,\ldots,k+1.$ By the operad
morphism property and the induction assumption we know that (\ref{P})
is true for each $S_j$.

\begin{lemma}
\label{case2}
The validity of  (\ref{P}) for the faces
$R_i$, $i=1,\ldots, k$,
implies its validity  for $T$.
\end{lemma}

\noindent
{\bf Proof.}
By definition of $R_i$ and $S_j$ and the property $\pa^2=0$,
\begin{eqnarray*}
\pa T^+ &=&   T + \sum   R_i + \sum   S_j, \\
0=\pa\pa T^+&=& \pa T + \sum \pa R_i + \sum \pa S_j.
\end{eqnarray*}
Therefore,
$$-\pa T = \sum   \pa R_i + \sum \pa S_j.$$
 Moreover,
$$p(\pa R_i)=\pa p( R_i)\quad \mbox{\rm and }\quad
p(\pa S_j)=\pa p( S_j).$$
Thus
\begin{eqnarray*}
p(-\pa T )&=& \sum   p(\pa R_i) + \sum   p(\pa S_j)
= \sum \pa p(R_i) + \sum   \pa p(S_j)=\pa p(\sum   R_i +\sum  S_j)\\
&=&\pa p(\pa T^+ - T)=\pa p(\pa T^+) -\pa p(T)
\pa\pa (T^+)-\pa p(T)=-\pa p(T).
\end{eqnarray*}
The summations in the above displays are taken
over $1 \leq i \leq k$ and $1 \leq j \leq k+1$.
This completes the proof of Lemma \ref{case2} and  the induction in
the proof of  (\ref{P}).\qed


\begin{thebibliography}{1}

\bibitem{boardman-vogt:73}
J.M. Boardman and R.M. Vogt.
\newblock {\em Homotopy Invariant Algebraic Structures on Topological Spaces}.
\newblock Springer-Verlag, 1973.

\bibitem{gaberdiel-zwiebach:NP97}
M.R. Gaberdiel and B.~Zwiebach.
\newblock Tensor constructions of open string theories {I}: Foundations.
\newblock {\em Nucl. Phys. B}, 505(3):569--624, Nov. 1997.

\bibitem{getzler-jones:preprint}
E.~Getzler and J.D.S. Jones.
\newblock Operads, homotopy algebra, and iterated integrals for double loop
  spaces.
\newblock Preprint {\tt hep-th/9403055}, March 1994.


\bibitem{ko-so}
M.~Kontsevich and Y.~Soibelman.
\newblock Homological mirror symmetry and torus fibrations.
\newblock Preprint {\tt math.SG/0011041}, November 2000.

\bibitem{markl:JPAA92}
M.~Markl.
\newblock A cohomology theory for {$A(m)$-algebras} and applications.
\newblock {\em J. Pure Appl. Algebra}, 83:141--175, 1992.

\bibitem{markl:zebrulka}
M.~Markl.
\newblock Models for operads.
\newblock {\em Comm. Algebra}, 24(4):1471--1500, 1996.

\bibitem{markl:ha}
M.~Markl.
\newblock Homotopy algebras are homotopy algebras.
\newblock Preprint {\tt math.AT/9907138}, to appear in Forum Matematicum, July
  1999.

\bibitem{markl-shnider-stasheff:book}
M.~Markl, S.~Shnider, and J.~D. Stasheff.
\newblock {\em Operads in Algebra, Topology and Physics}, volume~96 of {\em
  Mathematical Surveys and Monographs}.
\newblock American Mathematical Society, Providence, Rhode Island, 2002.

\bibitem{umble-saneblidze}
S.~Saneblidze and R.~Umble.
\newblock A diagonal on the associahedra.
\newblock Preprint {\tt math.AT/0011065}, November 2000.

\bibitem{serre:AM51}
J.-P. Serre.
\newblock Homologie singuli\'ere des espaces fibre\'es.
\newblock {\em Ann. of Math.}, 54:425--505, 1951.
\newblock In French.


\end{thebibliography}
\def\cprime{$'$}

\end{document}